\documentclass[12pt]{amsart}
\usepackage{amsmath,amssymb, amsthm,graphicx}

\usepackage{verbatim,enumerate}

\usepackage{url}
\usepackage{hyperref} 
\usepackage[capitalise]{cleveref} 

\graphicspath{{pictures/}{.0/pictures/}}
\usepackage{xcolor}
\usepackage{appendix}
\usepackage{mathtools}

\newcommand{\cvec}[1]{\begin{pmatrix}#1\end{pmatrix}}

\DeclarePairedDelimiter\ceil{\lceil}{\rceil}
\DeclarePairedDelimiter\floor{\lfloor}{\rfloor}
\DeclarePairedDelimiter\set{\{}{\}}

\DeclarePairedDelimiter\brac{[}{]}

\newcommand{\card}[1]{\ensuremath{\# #1}}  

\DeclareMathOperator{\conv}{conv} 

\newcommand{\R}{\mathbb{R}}

\newcommand{\f}{\varphi}

\newcommand{\ff}{\eta} 
\newcommand{\hh}{\tau} 
\newcommand{\ii}{\theta}
\newcommand{\jj}{\rho}

\renewcommand{\geq}{\geqslant} 
\renewcommand{\leq}{\leqslant} 
\renewcommand{\ge}{\geqslant} 
\renewcommand{\le}{\leqslant}

\newtheorem{theorem}{Theorem}
\newtheorem{corollary}[theorem]{Corollary}
\newtheorem{lemma}[theorem]{Lemma}
\newtheorem{proposition}[theorem]{Proposition}

\newtheorem{conjecture}[theorem]{Conjecture}
\crefname{conjecture}{Conjecture}{Conjectures}

\crefname{question}{Question}{Questions}

\theoremstyle{definition}

\newtheorem{fact}{Fact}
\AtEndEnvironment{proof}{\setcounter{fact}{0}}
\newenvironment{factproof}{\noindent\emph{Proof of fact.}}{\hfill$\qed$}

\theoremstyle{remark}
\newtheorem{remark}[theorem]{Remark}

\newtheorem{claim}{Claim}
\AtEndEnvironment{proof}{\setcounter{claim}{0}}
\newenvironment{claimproof}{\noindent\textsc{Proof of claim}}{\hfill$\qed$}
\crefname{claim}{Claim}{Claims}

\newtheorem{case}{Case}
\AtEndEnvironment{proof}{\setcounter{case}{0}}
\AtEndEnvironment{claimproof}{\setcounter{case}{0}}

\usepackage{anysize}
\marginsize{3cm}{3cm}{3cm}{3cm}

\begin{document}

\date{\today}
\title{A lower bound theorem for  $d$-polytopes with $2d+2$  vertices}
\author{Guillermo Pineda-Villavicencio}
\author{Aholiab Tritama}
\author{Jie Wang}

\address{School of Information Technology, Deakin University, Locked Bag 20000, Geelong VIC. 3220, Australia}
\email{\texttt{guillermo.pineda@deakin.edu.au},\texttt{a.tritama@research.deakin.edu.au}, \allowbreak
\texttt{wangjiebulingbuling@foxmail.com}}

\author{David Yost}
\address{Federation University, Mount Helen, Vic. 3350, Australia}
\email{\texttt{d.yost@federation.edu.au}}


\begin{abstract} We  establish a lower bound theorem for the number of $k$-faces ($1\le k\le d-2$) in a $d$-dimensional polytope $P$ (abbreviated as a $d$-polytope) with $2d+2$ vertices,  extending the previously known case for  $k=1$. We identify all minimisers for $d\le 5$.  Two distinct lower bounds emerge, depending on the number of facets of $P$. When $P$ has precisely $d+2$ facets, the lower bound is tight when $d$ is odd.  If  $P$ has at least $d+3$ facets, the lower bound is always tight, and  equality holds for some $1\le k\le d-2$  only when $P$ has precisely  $d+3$ facets. 

Moreover, for $1\le k\le \ceil{d/3}-2$, the minimisers among $d$-polytopes with $2d+2$ vertices have precisely $d+3$ facets, while for $\floor{0.4d}\le k\le d-1$, the lower bound arises from $d$-polytopes with $d+2$ facets. 
\end{abstract}

\maketitle
\markboth{\textsc{\small G.~Pineda-Villavicencio, A.~Tritama,  J.~Wang, and  D.~Yost}}{\shorttitle}
\section{Introduction}

Two (convex) polytopes  are combinatorially equivalent if their face lattices are isomorphic. We do not distinguish between combinatorially equivalent polytopes.

Gr\"unbaum~\cite[Sec.~10.2]{Gru03} conjectured that the function 
\begin{equation}\label{eq:function-at-most-2d}      
   \ii_k(d+s,d):=\binom{d+1}{k+1}+\binom{d}{k+1}-\binom{d+1-s}{k+1},\; \text{for  $1\le s\le d$}
\end{equation}
gives the minimum number of $k$-faces of a $d$-polytope with $d+s$ vertices and proved the conjecture for $s=2,3,4$. Xue~\cite{Xue21} then completed the proof, including a characterisation of the unique minimisers for $k\in[1\ldots d-2]$. Earlier, Pineda-Villavicencio et.~al~\cite{PinUgoYos15} had settled the conjecture for some values of $k$.

When $k\in [1\ldots d-2]$ (the set of integers between $1$ and $d-2$), each minimiser of \eqref{eq:function-at-most-2d} is a  $(d-s)$-fold pyramid over a  simplicial $s$-prism for $s\in [1\ldots d]$; we call such a polytope an \textit{($s$,$d-s$)-triplex} and denote it by $M(s,d-s)$.

Let $P$ be a $d$-polytope in $\R^{d}$, let $F$ be a face of $P$, and let $K$ be a closed halfspace in $\mathbb{R}^d$ such that the vertices of $P$ not in $K$ are the vertices of $F$. A polytope $P'$ is obtained by \textit{truncating the face} $F$ of $P$ if $P'=P\cap K$. 

For the parameters $d\ge 2$ and $s\in [2\ldots d]$, if we truncate a \textit{simple} vertex (a vertex contained in exactly $d$ facets) from the $(s,d-s)$-triplex, we obtain a $d$-polytope $J(s,d)$ with $2d+s-1$ vertices,  $d+3$ facets, and whose number of $k$-faces with $k\in [1\ldots d-1]$ is given by   the function
\begin{equation}
 \eta_k(2d+s-1,d):=\binom{d+1}{k+1}+2\binom{d}{k+1}-\binom{d+1-s}{k+1}.
\label{eq:3minus1-function}
\end{equation}
The paper~\cite{PinYos22} refers to $J(2,d)$ as a \textit{$d$-pentasm} and the polytope $J(3,d)$ as $B(d)$.

The {\it Cartesian product} of a $d$-polytope $P\subset \R^{d}$ and a $d'$-polytope $P'\subset \R^{d'}$ is the Cartesian product of the sets $P$ and $P'$: $P\times P'=\set{(p, p')^{t}\in \R^{d+d'}\mid p\in P,\, p'\in P'}$.  We denote by $T(d)$ the $d$-simplex and by $f_{k}$ the number of $k$-faces in a polytope. 

\begin{lemma}[McMullen 1971, {\cite[Sec.~3]{McM71}}]\label{lem:dplus2facets}
Let  $P$ be a $d$-dimensional polytope with $d+2$ facets, where $d\ge 2$. Then, there exist  integers $2\le a\le d$ and $1\le m\le \floor{a/2}$ such that $P$ is a $(d-a)$-fold pyramid over $T(m)\times T(a-m)$. The number of  $k$-faces of $P$ is
\begin{equation}\label{eq:dplus2facets}
\jj_k(a,m,d):= \binom{d+2}{k+2} -\binom{d-a+m+1}{k+2}-\binom{d-m+1}{k+2} +\binom{d-a+1}{k+2}.
\end{equation}
In particular,  $f_{0}(P)=d+1+m(a-m)$.
\end{lemma}
 
 Pineda-Villavicencio~\cite[Prob.~8.7.11]{Pin24} conjectured $J(s,d)$ and some $d$-polytopes with $d+2$ facets minimise the number of faces in polytopes with up to $3d-1$ vertices. Our results provide further evidence in support of this conjecture.
 
 \begin{conjecture}[{\cite[Prob.~8.7.11]{Pin24}}]
 \label{conj:3minus1}
   Let $d\ge 3$ and $s\in [2\ldots d]$. Let $P$ be a $d$-polytope with $2d+s-1$ vertices.
\begin{enumerate}[{\rm (i)}]
\item If $P$ has at least $d+3$ facets, then $f_{k}(P)\ge f_{k}(J(s,d))$ for each $k\in [1\ldots d-2]$.
\item  If $P$ has $d+2$ facets, then it is a $(d-a)$-fold pyramid over $T(m)\times T(a-m)$ for some $a\in [2\ldots d]$ and $2\le m\le \floor{a/2}$.
\end{enumerate}
 \end{conjecture}

The case $s=2$ in \cref{conj:3minus1} was settled independently by Pineda-Villavicencio and Yost~\cite{PinYos22}, and by Xue~\cite{Xue24}. Here we verify the case of $s=3$. Let 
\begin{equation*}
\begin{aligned}
    \hh_{k}(2d+2,d)&:=\binom{d + 1}{k + 1}+\binom{d}{k + 1}+\binom{d - 1}{k + 1}-\binom{\ceil{(d+1)/2}-1}{k + 1}\\
    &\quad -  \binom{\ceil{(d+1)/2}-2}{k + 1}.
\end{aligned}    
    \end{equation*}
See also \cref{lem:lower-bound-dplus2-facets}. Additionally, denote by $A(d)$ a polytope obtained by truncating a nonsimple vertex of  a $(2,d-2)$-triplex.   Our main theorem reads as follows.

\begin{theorem} Let $d\ge 3$, let $P$ be a $d$-polytope with $2d+2$ vertices, and let $k\in [1\ldots d-2]$. 

\begin{enumerate}[{\rm (i)}]
\item If $P$ has at least $d+3$ facets, then $f_{k}(P)\ge f_{k}(J(3,d))=f_{k}(A(d))=\eta_{k}(2d+2,d)$.  Moreover, equality holds for some $k$ only if $P$ has exactly $d+3$ facets.

\item If $P$ has $d+2$ facets, then  $f_{k}(P)\ge \hh_{k}(2d+2,d)$. Furthermore, the equality is tight for every $k\in [1\ldots d-2]$ when $d$ is odd, $a=\floor{(d+3)/2}+1$, and $P$ is a $(d-a)$-fold pyramid over $T(2)\times T(a-2)$. 
\end{enumerate}

\label{thm:2dplus2}
\end{theorem}

We remark that, when $d\ge 9$, another dichotomy manifests:   

\begin{enumerate}[{\rm (i)}] 
\item If $1\le k\le \ceil{d/3}-2$, then 
$\eta_{k}(2d+2,d)<\tau_{k}(2d+2,d)$ (\cref{lem:combinatorial-ineq}(i)), so $\eta_{k}(2d+2,d)$ is the lower bound for the number of $k$-faces in $d$-polytopes with $2d+2$ vertices. 

\item If $\floor{0.4d}\le k\le d-1$, then $\eta_{k}(2d+2,d)>\tau_{k}(2d+2,d)$ (\cref{lem:combinatorial-ineq}(ii)), making $\tau_{k}(2d+2,d)$ the lower bound for $k$-faces in $d$-polytopes with $2d+2$ vertices.
\end{enumerate}

It is worth noting that a $d$-polytope with $2d+2$ vertices and $d+2$ facets exists (if and) only if $d+1$ is composite and $d\ne3$. This is an easy consequence of the last sentence of \cref{lem:dplus2facets}.
 
The proof of Part (i) of \cref{thm:2dplus2} (restated as \cref{thm:2dplus2-restated}) is lengthy, involving case analysis, induction, properties of {\it super-Kirkman}  polytopes---where every pair of facets intersects in a ridge---and a key counting result of Xue~\cite[Prop.~3.1]{Xue21}, restated in \cref{cor:number-faces-outside-facet-practical},  which estimates the number of faces containing at least one vertex from a  set of at most $d$ vertices.

\section{Polytopes with $d+2$ vertices or facets}
\label{sec:dplus2} 
By duality, the following lemma for $d$-polytopes with $d+2$ vertices is equivalent to \cref{lem:dplus2facets}. Historically, \cref{lem:dplus2vertices} was proved first, and  \cref{lem:dplus2facets} deduced from it later.  The {\it direct sum} $P\oplus P'$ of a $d$-polytope $P\subset \R^{d}$ and a $d'$-polytope $P'\subset\R^{d'}$ with the origin in their relative interiors is the $(d+d')$-polytope: 
\begin{equation*}
\label{eq:direct-sums}
P\oplus P'=\conv \left(\left\{\cvec{ p\\0_{d'}}\in \R^{d+d'}\middle|\;  p\in P\right\}\bigcup \left\{\cvec{ 0_{d}\\ p'}\in \R^{d+d'}\middle|\;  p'\in P'\right\}\right).    
\end{equation*} 
Here, $0_{d}$ denotes the zero vector in $\R^{d}$.
\begin{lemma}\label{lem:dplus2vertices}\cite[Sec.~6.1]{Gru03}
Let  $P$ be a $d$-dimensional polytope with $d+2$ vertices, where $d\ge 2$. Then, there exist  integers $2\le a\le d$ and $1\le m\le \floor{a/2}$ such that $P$ is a $(d-a)$-fold pyramid over $T(m)\oplus T(a-m)$. The number of  $k$-dimensional faces of $P$ is
\begin{equation*}
\binom{d+2}{d-k+1} -\binom{d-a+m+1}{d-k+1}-\binom{d-m+1}{d-k+1} +\binom{d-a+1}{d-k+1}.
\end{equation*}
In particular,  $f_{d-1}(P)=d+1+m(a-m)$.
\end{lemma}

As in Gr{{\"u}}nbaum~\cite[Thm.~6.1.4]{Gru03}, we denote the $(d-a)$-fold pyramid over $T(m)\oplus T(a-m)$ by $T_{m}^{d,d-a}$. Gr{{\"u}}nbaum~\cite[p.~101]{Gru03} established inequalities for the number of $k$-faces in $d$-polytopes with $d+2$ vertices.

\begin{lemma} For $k\in [0\ldots d-1]$, the following hold:
\begin{enumerate}[{\rm (i)}]
\item If $2\le a\le d$ and $1\le m\le \floor{a/2}-1$, then  $f_{k}(T_{m}^{d,d-a})\le f_{k}(T_{m+1}^{d,d-a})$, with strict inequality if and only if $m\le k$.
\item If $2\le a\le d-1$ and $1\le m\le \floor{a/2}$, then  $f_{k}(T_{m}^{d,d-a})\le f_{k}(T_{m}^{d,d-(a+1)})$, with strict inequality if and only if $a-m\le k$.
\end{enumerate}
\label{lem:dplus2-vertices-inequalities}
\end{lemma}

The dual polytope $(T_{m}^{d,d-a})^{*}$ of $T^{d,d-a}_{m}$ is a $(d-a)$-fold pyramid over $T(m)\times T(a-m)$. We list the $d+2$ facets of $(T_{m}^{d,d-a})^{*}$ below.

\begin{remark}\label{rem:dplus2facets-facets}
For $d\ge 2$, $2\le a\le d$, and $1\le m\le \floor{a/2}$,  the $d+2$ facets of the $(d-a)$-fold pyramid $P$ over $T(m)\times T(a-m)$ are as follows:
\begin{enumerate}[{\rm (i)}]
\item $m+1$ facets that are  $(d-a)$-fold pyramids over $T(m-1)\times T(a-m)$, whose number of vertices is $f_{0}(P)-(a-m+1)$;

\item $a-m+1$ facets that are $(d-a)$-fold pyramids over $T(m)\times T(a-m-1)$, whose number of vertices is $f_{0}(P)-(m+1)$;

\item $d-a$ facets that are $(d-a-1)$-fold pyramids over $T(m)\times T(a-m)$, whose number of vertices is $f_{0}(P)-1$.
\end{enumerate}
\end{remark}

The next lemma, whose proof is embedded in the proof of \cite[Thm.~4.1]{Xue24}, settles \cref{thm:2dplus2}(ii). We write a proof for the sake of completeness. 
\begin{lemma} For $d\ge 4$, $s\in [2\ldots d-2]$,  $a=\floor{(d+s)/2}+1$,  and $1\le k\le d-2$, the number \begin{align*}
\hh_{k}(2d+s-1,d)&:=\binom{d+2}{k+2} -\binom{d-a+3}{k+2}-\binom{d-1}{k+2} +\binom{d-a+1}{k+2}\\
&\,=\binom{d+1}{k+1}+\binom{d}{k+1}+\binom{d-1}{k+1} -\binom{d-a+2}{k+1}-\binom{d-a+1}{k+1}
\end{align*}
of $k$-faces  of the $d$-polytope $(T_{2}^{d,d-a})^{*}$ is a lower bound on the number of $k$-faces of $d$-polytopes with exactly $d+2$ facets and \textbf{at least} $2d+s-1$ vertices. 
\label{lem:lower-bound-dplus2-facets}
\end{lemma}
	
\begin{proof} We reason as in \cite[Sec.~4]{Xue24}. Let $P:=(T_{2}^{d,d-a})^{*}$. Then $f_{0}(P)=d+2\floor{(d+s)/2}-1$. Now consider a $d$-polytope $Q$, distinct from $P$, with $d+2 $ facets and at least $2d+s-1$ vertices. By \cref{lem:dplus2facets},
$Q$ must have the form $Q=(T_{m}^{d,d-b})^{*}$,  where $2\le b\le d$ and $1\le m\le \floor{b/2}$. Since $f_0(Q)=d+1+m(b-m)\ge 2d+s-1$ and $b\le d$, we must have $m\ge 2$.

Assume that $f_{0}(Q)\ge f_{0}(P)$. We want to prove that $f_{k}(P)\le f_{k}(Q)$ for $1\le k\le d-2$.  We work in the dual setting, proving instead that $f_{k}(P^{*})\le f_{k}(Q^{*})$ for $1\le k\le d-2$. From $f_{d-1}(Q^{*})\ge f_{d-1}(P^{*})$, we get
\begin{equation}
m(b-m)\ge 2(a-2).
\label{eq:lower-bound-dplus2-facets-1}
\end{equation}
Suppose that $m=2$. Then $b\ge a$ and \cref{lem:dplus2-vertices-inequalities}(ii) yields that $f_{k}(P^{*})\le f_{k}(Q^{*})$ for $1\le k\le d-2$, as desired. Consequently, we assume that $m\ge 3$. Hence $b\ge 6$.

If $b\ge a$ then  \cref{lem:dplus2-vertices-inequalities}(ii) again yields that $f_{k}(T_{m}^{d,d-a})\le f_{k}(Q^{*})$ for $1\le k\le d-2$. Additionally, $f_{k}(T_{m}^{d,d-a})\ge f_{k}(P^{*})$ for $1\le k\le d-2$ by \cref{lem:dplus2-vertices-inequalities}(i). Hence, we are home in this case. Consequently, we may assume that $a> b$.
\begin{equation}
\label{eq:lower-bound-dplus2-facets-1a}
\begin{aligned}
f_{k}(Q^{*})-f_{k}(P^{*})&=\underbrace{\brac*{-\binom{d-b+m+1}{d-k+1}+\binom{d-a+3}{d-k+1}}}_{=A}\\
&+\underbrace{\brac*{-\binom{d-m+1}{d-k+1}+\binom{d-1}{d-k+1}}}_{=B} \\
&+\underbrace{\brac*{\binom{d-b+1}{d-k+1}-\binom{d-a+1}{d-k+1}}}_{=C}.
\end{aligned}
\end{equation}
Using \cref{lem:combinatorial-identities}(iii), we write $A,B,C$ as follows:
\begin{gather*}
A=-\sum_{\ell=1}^{(a-b)+m-2}\binom{d-b+m+1-\ell}{d-k},\; B=\sum_{\ell=1}^{m-2}\binom{d-1-\ell}{d-k},\\  C=\sum_{\ell=1}^{a-b}\binom{d-b+1-\ell}{d-k}.
\end{gather*}
From this we get expressions for $A+C$ and $A+B+C$.
\begin{align*}
A+C&=-\sum_{\ell=1}^{(a-b)+m-2}\binom{d-b+m+1-\ell}{d-k}+\sum_{\ell=1}^{a-b}\binom{d-b+1-\ell}{d-k}\\
&= -\sum_{\ell=1}^{m}\binom{d-b+m+1-\ell}{d-k}+\sum_{\ell=1}^{2}\binom{d-a+\ell}{d-k}\\
&\quad +\underbrace{\brac*{-\sum_{\ell=m+1}^{a-b+m-2}\binom{d-b+m+1-\ell}{d-k}+\sum_{\ell=1}^{a-b-2}\binom{d-b+1-\ell}{d-k}}}_{=0}\\
&=-\sum_{\ell=1}^{m}\binom{d-b+\ell}{d-k}+\sum_{\ell=1}^{2}\binom{d-a+\ell}{d-k}\\
&=-\sum_{\ell=1}^{m-2}\binom{d-b+\ell+2}{d-k}-\sum_{\ell=1}^{2}\brac*{\binom{d-b+\ell}{d-k}-\binom{d-a+\ell}{d-k}}.
\end{align*}
\begin{align*}
A+C+B&=-\sum_{\ell=1}^{m-2}\binom{d-b+\ell+2}{d-k}-\sum_{\ell=1}^{2}\brac*{\binom{d-b+\ell}{d-k}-\binom{d-a+\ell}{d-k}}\\
&\quad +\sum_{\ell=1}^{m-2}\binom{d-1-\ell}{d-k}\\
&=\sum_{\ell=1}^{m-2}\brac*{\binom{d-1-\ell}{d-k}-\binom{d-b+\ell+2}{d-k}}-\sum_{\ell=1}^{2}\brac*{\binom{d-b+\ell}{d-k}-\binom{d-a+\ell}{d-k}}.  
\end{align*}	
Using \cref{lem:combinatorial-identities}(iii), we further get 
\begin{align*}
A+C+B=\underbrace{\sum_{\ell=1}^{m-2}\brac*{\sum_{h=1}^{b-3-2\ell}\binom{d-1-\ell-h}{d-k-1}}}_{=D}-\underbrace{\sum_{\ell=1}^{2}\brac*{\sum_{h=1}^{a-b}\binom{d-b+\ell-h}{d-k-1}}}_{=E}.
\end{align*}
The coefficient $\binom{d-b+3}{d-k-1}$ is the smallest among those in $D$. Thus, 
\begin{equation}
\label{eq:lower-bound-dplus2-facets-2}
D\ge \sum_{\ell=1}^{m-2} (b-3-2\ell)\binom{d-b+3}{d-k-1}.
\end{equation}
The coefficient $\binom{d-b+1}{d-k-1}$ is the largest among those in $E$. Thus, 
\begin{equation}
\label{eq:lower-bound-dplus2-facets-3}
E\le 2(a-b)\binom{d-b+1}{d-k-1}.
\end{equation}
Combining \eqref{eq:lower-bound-dplus2-facets-2} and \eqref{eq:lower-bound-dplus2-facets-3}, we get that 
\begin{align*}
A+B+C&\ge \sum_{\ell=1}^{m-2} (b-3-2\ell)\binom{d-b+3}{d-k-1}-2(a-b)\binom{d-b+1}{d-k-1}\\
&= (m-2)(b-2-m)\binom{d-b+3}{d-k-1}-2(a-b)\binom{d-b+1}{d-k-1}.
\end{align*}
Furthermore, \eqref{eq:lower-bound-dplus2-facets-1} is equivalent to $(m-2)(b-2-m)\ge 2(a-b)$, which implies that 
\begin{equation}
\label{eq:lower-bound-dplus2-facets-4}
A+B+C\ge 2(a-b)\brac*{\binom{d-b+3}{d-k-1}-\binom{d-b+1}{d-k-1}}\ge 0,\ \text{(since $a>b$)}.
\end{equation}
The lemma follows.
\end{proof}
The function $\tau_{k}(2d+s-1,d)$ in \cref{lem:lower-bound-dplus2-facets} is obtained by setting $a=\floor{(d+s)/2}+1$ and $m=2$ in the function $\rho_{k} (a,m,d)$ from \cref{lem:dplus2facets}. 

 


\section{Known lower bounds  for polytopes with few vertices and auxiliary results}

We define a vertex figure of a polytope $P$ at a vertex $v$ as the polytope $P/v:=H\cap P$ where $H$ is a hyperplane  separating $v$ from the other vertices of $P$. There is a bijection between the $k$-faces of $P$ that contain $v$ and the $(k-1)$-faces of $P/v$.  We require a result by Xue~\cite[Prop.~3.1]{Xue21}; see also  \cite[Prop.~12]{PinYos22}. 

\begin{proposition}[Xue 2021]\label{prop:number-faces-outside-facet}
Let $d\ge 2$ and let $P$ be a $d$-polytope. In addition, suppose that  $r\le d+1$ is given and  that $S:=(v_1,v_2,\ldots,v_r)$ is a sequence of distinct vertices in $P$. Then the following hold.
\begin{enumerate}[{\rm (i)}]
    \item  There is a sequence $F_1, F_2,\ldots, F_r$ of faces of $P$ such that each $F_i$ has dimension $d-i+1$ and contains $v_i$, but does not contain any $v_j$ with $j<i$.
    \item For each $k\ge1$, the number of $k$-faces of $P$ that contain at least one of the vertices in $S$ is bounded from below by \[\sum_{i=1}^{r}f_{k-1}(F_i/v_i).\]
\end{enumerate}
\end{proposition}

	Xue's lower bound theorem for $d$-polytopes with at most $2d$ vertices (see \cite{Xue21}) reads as follows:
  
\begin{theorem}[$d$-polytopes with at most $2d$ vertices,  {\cite{Xue21}}]\label{thm:at-most-2d}  
Let $d\ge 2$ and $1\le s\le d$. If $P$ is a  $d$-polytope with  $d+s$ vertices, then
\[f_k(P)\ge \ii_k(d+s,d),\;  \text{for all $k\in[1\ldots d-1]$}.\] 
Also, if $f_k(P)=\ii_k(d+s,d)$ for some $k\in[1\ldots d-2]$, then $P$ is the $(s,d-s)$-triplex.
\end{theorem}

We proceed with a lower bound theorem for $d$-polytopes with $2d+1$ vertices, independently proved  by Pineda-Villavicencio and Yost~\cite{PinYos22}, and by Xue~\cite{Xue24}. Denote by $f(P)$ the $f$-vector of a polytope $P$, and denote by $\Sigma(d)$ a $d$-polytope combinatorially equivalent to the convex hull of
 \[\set{0,e_1,e_1+e_k,e_2,e_2+e_k,e_1+e_2,e_1+e_2+2e_k: 3\le k\le d},\]
where $\{e_i\}$ is the standard basis of $\mathbb{R}^d$. 

\begin{theorem}[Pineda-Villavicencio and Yost 2022; Xue 2023]
\label{thm:2dplus1}
Let $d\ge 3$, $P$  a $d$-polytope with \textbf{at least} $2d+1$ vertices, and $k\in [1\ldots d-2]$. 
\begin{enumerate}[{\rm (i)}]
\item Let $d=3$. If $P$ is $\Sigma(3)$ or a 3-pentasm, then $f(P)=(7,11,6)$. Otherwise, $f_1(P)>11$ and $f_2(P)>6$.
\item  If $d\ge 4$ and $P$ has at least $d+3$ facets, then $f_{k}(P)\ge f_{k}(J(2,d))$, with equality for some $k\in [1\ldots d-2]$ only if  $P=J(2,d)$.
\item  If $d\ge 4$ and $P$ has  $d+2$ facets, then $f_{k}(P)\ge f_{k}((T_{2}^{d,d-(\floor{d/2}+2)})^{*})$. If $d$ is even  then $(T_{2}^{d,d-(\floor{d/2}+2)})^{*}$ has $2d+1$ vertices.
\end{enumerate}
\end{theorem}

A corollary of \cref{thm:at-most-2d,thm:2dplus1} extends to $d$-polytopes with at least $d+s$ vertices. 

\begin{corollary} For each $d\ge2$, each $k\in [1\ldots d-1]$, and each $s\in[1\ldots d]$, if $P$ is a $d$-polytope with \textbf{at least} $d+s$ vertices, then $f_{k}(P)\ge \theta_{k}(d+s,d)$.

Moreover, if $f_k(P)=\ii_k(d+s,d)$ for some $k\in[1\ldots d-2]$, then $P$ has $d+2$ facets.
\label{cor:more-dpluss}
\end{corollary}
\begin{proof}
 The case  $d=2$ is trivial. Let $d\ge 3$. If $P$ has at most $2d$ vertices, the claim follows from \cref{thm:at-most-2d}. Hence assume $f_0(P)\ge 2d+1$. Since $\ii_k(d+s, d)$ is an increasing function in $s$, it suffices to show that $f_k(P) \ge \ii_k(2d, d)=\binom{d+1}{k+1}+\binom{d}{k+1}$ for each $k\in [1\ldots d-2]$.
 If $P$ has at least $d+3$ facets, then by \cref{thm:2dplus1}(ii) we have
 \begin{equation*}
     f_k(P)\ge \eta_k(2d+1,d)=\binom{d+1}{k+1}+\binom{d}{k+1}+\binom{d-1}{k}>\ii_k(2d, d).
 \end{equation*}
 If $P$ has exactly $d+2$ facets, then, by \cref{thm:2dplus1}(iii) (see also \cref{lem:lower-bound-dplus2-facets}),  
 \begin{align*}
f_k(P)&\ge f_{k}((T_{2}^{d,d-(\floor{d/2}+2)})^{*}) = \binom{d+1}{k+1}+\binom{d}{k+1}+\binom{d-1}{k+1}-\binom{\ceil{d/2}}{k+1}\\
&\quad-\binom{\ceil{d/2}-1}{k+1} \\
&> \binom{d+1}{k+1}+\binom{d}{k+1} = \ii_k(2d, d).
 \end{align*}
If $d$ is even, note that $d-1=\ceil{d/2}+\ceil{d/2}-1$, so the strict inequality is rather straightforward. If $d$ is odd, the polytope $(T_{2}^{d,d-(\floor{d/2}+2)})^{*}$ has exactly $2d$ vertices, and the strict inequality follows from \cref{thm:at-most-2d}. 
 
Additionally, if $f_k(P)=\ii_k(d+s,d)$, then by \cref{thm:at-most-2d}, $P$ must be a triplex with at least $d+s$ vertices, which has $d+2$ facets. This completes the proof.
\end{proof}

It is convenient to write a corollary of \cref{prop:number-faces-outside-facet} and \cref{cor:more-dpluss}.

\begin{corollary}\label{cor:number-faces-outside-facet-practical}
Let $d\ge 2$ and let $P$ be a $d$-polytope. In addition, suppose that  $r\le d+1$ is given and  that $S:=(v_1,v_2,\ldots,v_r)$ is a sequence of distinct vertices in $P$.  Then the following hold.
\begin{enumerate}[{\rm (i)}]
    \item  There is a sequence $F_1, F_2,\ldots, F_r$ of faces of $P$ such that each $F_i$ has dimension $d-i+1$ and contains $v_i$, but does not contain any $v_j$ with $j<i$.
    \item For each $1\le i\le r$, let $s_i$  satisfy $1\le s_{i} \le d-i+1$, and suppose that $\deg_{F_{i}}(v_i)\ge d-i+1+s_{i}$. Then, for each $k\ge1$, the number of $k$-faces of $P$ that contain at least one of the vertices in $S$ is bounded below by \[\sum_{i=1}^{r}f_{k-1}(F_i/v_i)\ge \sum_{i=1}^{r}\theta_{k-1}(d-i+1+s_{i},d-i).\]
\end{enumerate}
\end{corollary}

\section{Polytopes with $2d+2$ vertices: Small cases and pyramids over simple polytopes}

Denote by $C(d)$ a polytope obtained by truncating a {\it simple edge}, an edge with two simple vertices, of a  $(2,d-2)$-triplex. It has $3d-2$ vertices and $d+3$ facets.  Obviously $C(2)$ is just another quadrilateral.

\begin{theorem}[Small cases]\label{thm:2d+2-small} For $d=3,4,5$, the $d$-polytopes with exactly $2d+2$ vertices that minimise the number of $k$-faces for each $k\in [1\ldots d-1]$ are as follows.
\begin{enumerate}[{\rm (i)}]
\item For $d=3$, the minimisers are the polytopes $A(3)$ (the cube) and $J(3,3)$,  each with precisely $\eta_k(8,3)$ $k$-faces for $k=1,2$.
\item For $d=4$, the minimisers are the four polytopes $A(4)$, $J(3,4)$, $C(4)$, and $\Sigma(4)$, all with $f$-vector $(10, 21, 18, 7)$; equivalently, each has precisely $\eta_k(10,4)$ $k$-faces for $k=1,2,3$.
\item For $d=5$, the minimiser with seven facets is $T(2)\times T(3)$, with $f$-vector $(12, 30, 34, 21, 7)$. 

Among the 5-polytopes with at least eight facets, the minimisers include $A(5)$ and $J(3,5)$, both with $f$-vector $(12, 32, 39, 25, 8)$; equivalently, each has precisely $\eta_k(12,5)$ $k$-faces for $k=1,2,3,4$. These two are the only minimisers for $k=1,2$, but not for $k=3,4$. 
\end{enumerate}
Moreover, in each case, for every $k\in [1\ldots d-1]$, the minimisers  have  either $d+2$ or $d+3$ facets.
 \end{theorem}
\begin{proof}
(i)--(ii) The case $d=3$ of \cref{thm:2d+2-small} is  easy to check; one may also consult the catalogue \cite{BriDun73}.   The 4-polytopes with 10 vertices and the minimum number of edges, namely 21, are the polytopes $A(4)$, $J(3,4)$, $C(4)$, and $\Sigma(4)$ \cite[Thm.~6.1]{PinUgoYos16a}.
The Euler--Poincar\'e-Schl\"afli equation for $d=4$ \cite[Thm.~2.12.17]{Pin24} yields that  
\begin{equation*}
f_{0}-f_{1}+f_{2}-f_{3}=0. 
\end{equation*}
If $f_{0}=10$ and $f_{1}\ge 22$, this equation implies  that $f_{2}-f_{3}\ge 12$, or equivalently  that $f_{2}\ge 12+f_{3}$. Since  no $4$-polytope with 10 vertices has $6$ facets (\cref{lem:dplus2facets}), the inequality $f_{3}\ge 7$ gives that $f_{2}\ge 19$. This settles the case $d=4$. 

(iii) Suppose that $d=5$. The unique $5$-polytope with 7 facets and 12 vertices is $T(2)\times T(3)$ (\cref{lem:dplus2facets}).   Consider a 5-polytope $P$ with $12$ vertices and at least 8 facets. We show the following.
\begin{enumerate}[{\rm (a)}]
\item  $f_1(P)\ge32$, with equality only for $A(5)$ and $J(3,5)$;

\item  $f_2(P)\ge39$, with equality only for $A(5)$ and $J(3,5)$;

\item  $f_3(P)\ge25$ but there are other minimisers.
\end{enumerate}

(a) According to \cite[Thm.~13]{PinUgoYos22}, the 5-polytopes with 12 vertices, at least 8 facets, and the minimum number of edges, namely 32, are the polytopes $A(5)$ and $J(3,5)$. 	

(c) We consider the dual setting: if $P$ has 12 vertices and at least 8 facets, then the dual polytope $P'$ has 12 facets and at least 8 vertices. Suppose first that $f_{0}(P')=8$. According to \cref{thm:at-most-2d}, the unique  $5$-polytope with $8$ vertices and at most $\theta_{1}(8,5)=22$ edges is the $(3,2)$-triplex, which has $7$ facets.  Furthermore, by reasoning with Gale transforms \cite[Sec.~2.14]{Pin24} or consulting the catalogues in \cite{FukMiyMor13a}, we can deduce  that the only 5-polytopes with 8 vertices and exactly $\theta_{1} (8, 5) + 1=23$ edges are the 3-fold pyramid over the pentagon and the 2-fold pyramid over the tetragonal antiwedge, each with 8 facets. The \textit{tetragonal antiwedge} is the unique nonpyramidal 3-polytope with six vertices and six 2-faces. The $d$-polytopes with $d+3$ vertices and $\theta_{1}(d+3,d)+2$ edges are characterised  in \cite[Prop.~16 (i)--(vi)]{PinYos22}; they all have $2d+1$ or fewer facets; for the particular case of $d=5$, we can also consult the catalogues in \cite{FukMiyMor13a}. Dually, we see that a 5-polytope with 8 facets and no more than $\theta_{1}(8,5)+2=24$ ridges has at most 11 vertices. Thus, $f_{3}(P) \ge 25=f_3(A({5}))=f_{3}(J(3,5))$, and also $f_{3}(P)-f_{4}(P)\ge 25-8=17$.

Suppose now that $P$ has 9 facets. Then, combining  \cref{thm:at-most-2d} and \cite[Thm.~20]{PinUgoYos15} we get that $P'$ has at least $\theta_{1}(9,5)+2=26$ edges. If instead $P$ has 10 facets, then, combining \cref{thm:at-most-2d} and \cite[Prop.~17]{PinUgoYos15} we get that the dual $P'$ has at least $\theta_{1}(10,5)+2=27$ edges. Hence, $f_{3}(P)\ge 26$ and $f_{3}(P)-f_{4}(P)\ge 17$ in both cases. Finally, suppose that $P$ has at least $11$ facets. In this case, the dual graph of $P$ informs us that  $2f_3(P)\ge 5f_4(P)\ge55$; that is, $f_3\ge28$ and  $f_3-f_4\ge1.5f_4\ge16.5$. From this analysis, regardless of the number of facets in $P$, we get that $f_3(P)\ge25$ again. Note for subsequent use that 
\begin{equation}
\label{eq:2d+2-small-1}
f_{3}(P)-f_{4}(P)\ge 17.
\end{equation}
in each case. 

There are three other minimisers of $f_3$ \cite{PinWanYos25}. The easiest one to describe is a pyramid over $J(4,4)$ with $f=(12, 33, 40, 25, 8)$.

(b) If $P$ is not $A(5)$ or $J(3,5)$, then $f_{1}(P)\ge 33$ from (a), in which case \eqref{eq:2d+2-small-1} and the Euler--Poincar\'e-Schl\"afli equation yields that  
\begin{equation*}
f_2=f_1-f_0+f_3-f_4+2\ge33-12+17+2=40.
\end{equation*}
This concludes the proof of the theorem. 
\end{proof}

\subsection{Pyramids with at least $d+3$ facets and simple base polytope}

We require the lower bound theorem for simple polytopes~\cite{Bar71,Bar73} and the notion of \textit{truncation polytopes},  polytopes obtained from simplices by successive truncation of vertices. 

\begin{remark}  The smallest vertex counts of truncation $d$-polytopes are $d + 1$, $2d$, $3d-1$, and $4d-2$.
\label{rmk:truncation-polytopes-vertices}
\end{remark}

\begin{theorem}[Simple polytopes, Barnette (1971--73)]
\label{thm:simple-lbt} 
Let $d\ge 2$ and let $P$ be a simple $d$-polytope with $f_{d-1}$ facets. Then
\begin{equation}\label{eq:simple}
 f_k(P)\ge\begin{cases}
(d-1)f_{d-1}-(d+1)(d-2),& \text{if $k=0$}; \\
\binom{d}{k+1}f_{d-1}-\binom{d+1}{k+1}(d-1-k),& \text{if $k\in [1\ldots d-2]$}.
\end{cases}
\end{equation}
If, for $d\ge 4$,  $f_{k}(P)$ achieves equality for some $k\in [0\ldots d-2]$, then $P$ must be a truncation polytope. For $d=2,3$, equality holds for every simple $d$-polytope.
\end{theorem}

 \begin{lemma} Let $d\ge 3$ and $0\le t\le d-2$. Additionally, let $P$ be a $t$-fold pyramid over a simple $(d-t)$-polytope and suppose $P$ has \textbf{at least} $2d+2$ vertices and \textbf{at least}  $d+3$ facets. Then, for  each $k\in [1\ldots d-2]$, the following hold:  \begin{enumerate}[{\rm (i)}]
 \item If $t=0$, then $f_{0}(P)\ge 3d-1$ and $f_{k}(P)\ge \eta_{k}(3d-1,d)$. Equality  holds for some $k\in [1\ldots d-2]$ if and only if $P$ is either a cube or the simple $d$-polytope $J(d,d)$ ($d\ge 3$).

 \item If $t\ge 1$, then $f_{k}(P) \ge \eta_{k}(2d+2,d)$. Moreover equality holds if and only if $k=d-2$,  and $P$ is a $t$-fold pyramid over the simple $(d-t)$-polytope $J(d-t,d-t)$ with $t$ satisfying  $f_0(P)=3d-2t-1\ge 2d+2$.
 \end{enumerate}
 \label{lem:2d+s-dplus3-simple-pyramid}
 \end{lemma}
\begin{proof} 

(i) First suppose that $t=0$. Since $f_{d-1}(P)\ge d+3$, the lower bound theorem for simple polytopes yields that 
$f_{0}(P)\ge 3d-1$ and, for  $k\ge 1$, that 
  
\begin{align}
\label{eq:2d+s-dplus3-simple-pyramid-1}
f_{k}(P)&\ge \binom{d}{k+1}(d+3)-\binom{d+1}{k+1}(d-1-k)\\
&=2\binom{d}{k+1}+\binom{d+1}{k+1}+\underbrace{\brac*{(d+1)\binom{d}{k+1}-\binom{d+1}{k+1}(d-k)}}_{=0}\nonumber\\
&=2\binom{d}{k+1}+\binom{d+1}{k+1}=\eta_k(3d-1,d)+\binom{1}{k+1}=\eta_k(3d-1,d)\nonumber.
\end{align}
If $d=3$, $f_1(P)=\eta_1(3d-1,d)$ for all simple 3-polytopes with $3d-1$ vertices and $d+3$ faces; the cube and $J(3,3)$ are the only such polytopes \cite[Lem.~2.19(v)]{PinUgoYos16a}. Suppose $d\ge 4$. By \cite[Lem.~2.19(v)]{PinUgoYos16a}, $J(d,d)$ is the only simple $d$-polytope with $3d-1$ vertices and $d+3$ facets. If $f_{0}(P)> 3d-1$, then, by the Lower bound theorem for simple polytopes,  \eqref{eq:2d+s-dplus3-simple-pyramid-1} is a strict inequality, since a truncation $d$-polytope with  $d+3$ facets has $3d-1$ vertices; see also  \cref{rmk:truncation-polytopes-vertices}.

(ii) We first deal with $d=3$. If $t=2,3$,  then $P$ is a 3-simplex, which has fewer than 8 vertices. If $t=1$, then $P$ is a pyramid over an $n$-gon with $n\ge 7$, in which case we have that $f_{1}(P)=2n>12=\eta_{1}(8,3)$ and $f_{2}(P)=n+1>d+3$. Hence the case $d=3$ is settled.

Suppose $d=4$. We reason as in the case $d=3$. If $t=3,4$,  then $P$ is a 4-simplex, which has fewer than 10 vertices. If instead $t=2$, then $P$ is a two-fold pyramid over an $n$-gon with $n\ge 8$, in which case we have that $f_{1}(P)=3n+1>21=\eta_{1}(10,4)$ and $f_{2}(P)=3n+1>18=\eta_{2}(10,4)$. If $t=1$, then $P$ is a pyramid over a simple 3-polytope $F$ with at least ten vertices and at least six faces (there are no simple 3-polytopes with nine vertices). By Euler's polyhedral formula, $f_1(F)\ge 14$. This implies that $f_k(P)>\eta_{k}(10,4)$ for $k=1,2$. We remark that if $P$ were a pyramid over a cube or $J(3,3)$, it would have nine vertices, contrary to our assumption on the number of vertices.  Thus we can start our induction on $d$ for all $t\ge 0$. Assume $d\ge 5$.
\begin{claim} If $Q$ is a $t$-fold pyramid over  $J(d-t,d-t)$, then  
\[f_k(Q)=
\binom{d+1}{k+1}+2\binom{d}{k+1}-2\binom{t+1}{k+1}.
\]
(Of course the latter term is 0 when $k>t$.)
Furthermore, if  $f_0(Q)=3d-2t-1\ge 2d+2$, then
$f_{d-2}(Q)=\eta_{d-2}(2d+2,d)$, while $f_{k}(Q)>\eta_{k}(2d+2,d)$ for all $1\le k\le d-3$.
    \label{cl:faces-multifold-pyramid}
\end{claim}

\begin{claimproof} The face numbers of $Q$ follow directly from the multifold-pyramid formula; see, for instance, \cite[Thm.~4.2.2]{Gru03}. If  $k>t$, $f_{d-2}(Q)=\eta_{d-2}(2d+2,d)$; moreover,   $f_{k}(Q)>\eta_{k}(2d+2,d)$ for each $1\le k\le d-3$. Suppose $1\le k\le t$.   Since $2t+1\le d-2$, we must have  $1\le k<d-3$, in which case
 \begin{equation*}
     f_k(Q)-\eta_k(2d+2,d)=\binom{d-2}{k+1}-2\binom{t+1}{k+1}\ge \binom{2t+1}{k+1}-2\binom{t+1}{k+1}>0.
 \end{equation*}
Applying the Vandermonde's identity to $\binom{2t+1}{k+1}$ (\cref{lem:combinatorial-identities}(v)) gives the last strict inequality.
\end{claimproof}

Since $t\ge 1$, the polytope $P$ is a pyramid over a $(d-1)$-polytope $F$, itself a $(t-1)$-fold pyramid with at least $2(d-1)+2$ vertices and at least $(d-1)+3$ facets. Thus our induction hypothesis applies to $F$. It also follows that 
\begin{equation*}
\label{eq:2d+s-dplus3-simple-pyramid-2}
f_{k}(P)=f_{k}(F)+f_{k-1}(F).
\end{equation*}

\textbf{Suppose that $t\ge 2$}. By the induction hypothesis, we get that
\begin{equation}
\label{eq:2d+s-dplus3-simple-pyramid-3}
f_{k}(F)\ge \eta_{k}(2(d-1)+2,d-1)=\binom{d}{k+1}+2\binom{d-1}{k+1}-\binom{d-3}{k+1},
\end{equation}
with equality if and only if $k=d-3$, $F$ is a $(t-1)$-fold pyramid over $J(d-t,d-t)$, and $f_0(F)\ge 2(d-1)+3>2(d-1)+2$.

Since $f_{0}(F)=f_{0}(P)-1$, $f_{0}(P)\ge 2d+2$, and $d\ge 5$, the induction hypothesis  gives 
\begin{align*}
f_{1}(P)&>\eta_{1}(2(d-1)+2,d-1)+2d+1\\
&=\brac*{\binom{d}{2}+2\binom{d-1}{2}-\binom{d-3}{2}}+\brac*{\binom{d}{1}+2\binom{d-1}{1}-\binom{d-3}{1}}\\
&=\eta_{1}(2d+2,d).
\end{align*}
Furthermore, for $k\in [2\ldots d-2]$, we have 	
\begin{align*}
f_{k}(P)&\ge \eta_{k}(2(d-1)+2,d-1)+\eta_{k-1}(2(d-1)+2,d-1)=\eta_{k}(2d+2,d).
\end{align*}		 
For $2\le k\le d-3$, $f_{k-1}(F)>\eta_{k-1}(2(d-1)+2,d-1)$ by the induction hypothesis, whence $f_{k}(P)>\eta_{k}(2d+2,d)$. 

If $k=d-2$, the equality $f_{d-2}(P)=\eta_{d-2}(2d+2,d)$ would then imply that $f_{d-3}(F)=\eta_{d-3}(2(d-1)+2,d-1)$, which, by the induction hypothesis, holds if and only if $F$ is a $(t-1)$-fold pyramid over $J(d-t,d-t)$ and $f_0(F)\ge 2d$; observe that $f_{d-2}(F)=d+2=\eta_{d-2}(2(d-1)+2,d-1)$. Therefore, $f_{d-2}(P)=\eta_{d-2}(2d+2,d)$ only when $P$ is a $t$-fold pyramid over $J(d-t,d-t)$.

\textbf{Assume that $t=1$}. Then $F$ is a simple $(d-1)$-polytope, in which case (i) yields that $f_{0}(F)\ge 3(d-1)-1$ and $f_{0}(P)\ge 3d-3$. 
Moreover, as $f_{0}(F)\ge 3(d-1)-1$, Part (i) yields that $f_{k}(F)\ge \eta_{k}(3d-4,d-1)$ for each $k\in [1\ldots d-2]$. Therefore, for $2\le k\le d-2$, we find that 
\begin{align*}
f_{k}(P)&=f_{k}(F)+f_{k-1}(F)\\
&\ge \eta_{k}(3(d-1)-1,d-1)+\eta_{k-1}(3(d-1)-1,d-1)\\
& = \binom{d+1}{k+1}+2\binom{d}{k+1}-\binom{2}{k+1}
\ge \binom{d+1}{k+1}+2\binom{d}{k+1}-\binom{d-2}{k+1}\\
&= \eta_{k}(2d+2,d).
\end{align*}	
The last inequality is tight only for $k = d - 2$. By Part (i), if $d \ge 5$, then the first inequality is tight only when $F=J(d-1,d-1)$. Hence in the equality case,  $P$ is a pyramid over $J(d-1,d-1)$. Additionally, 
\begin{align*}
f_{1}(P)& \ge f_{1}(F)+3d-4\ge \binom{d}{2}+2\binom{d-1}{2}+\binom{d}{1}+2\binom{d-1}{1}-2\\
&=\eta_{1}(3d-3,d)+1>\eta_{1}(2d+2,d).
\end{align*}
This concludes the proof of the lemma.	
\end{proof}

There is a basic bound for the number of $k$-faces in a polytope.

\begin{lemma} If $F$ is a facet of a $d$-polytope $P$, then, for each $k\in [0\ldots d-1]$,  $f_{k}(P)\ge f_{k}(F)+f_{k-1}(F)$,  
\label{lem:facet-k-faces}
\end{lemma}
\begin{proof} For $k\ge 0$, each $(k-1)$-face $F_{0}$ in $F$  is contained in a  $k$-face $F_{1}$ not in $F$ such that $F_{0}=F_{1}\cap F$. The result now follows.
\end{proof}
	
\section{Proof of the main theorem}
\label{sec:main-proof}

We restate Part (i) of \cref{thm:2dplus2} for convenience.

\begin{theorem}[Restatement of {\cref{thm:2dplus2}(\rm{i})}]Let $d\ge 3$,  and let $P$ be a $d$-polytope with $2d+2$ vertices. If $P$ has at least $d+3$ facets, then $f_{k}(P)\ge f_{k}(J(3,d))=f_{k}(A(d))=\eta_{k}(2d+2,d)$, for each $k\in [1\ldots d-2]$. Moreover, equality holds for some $k$ only if $P$ has exactly $d+3$ facets.
 \label{thm:2dplus2-restated}
\end{theorem}

We prove the inequality first  and address the equality afterwards. The cases $d=3,4,5$ are settled in \cref{thm:2d+2-small}. Assume that $d\ge 6$. The paper \cite{PinUgoYos22} settled the case $k=1$ for all $d$.  Our proof of the inequality part proceeds  by induction on $d$ for all $k\in [2\ldots d-2]$.  

\subsection{Claim \ref{cl:2d+2-dplus3-1}: A facet with $f_{0}(P)-2$ vertices}

We begin with an overview of the proof. Let $F$ be a facet with $f_{0}(P)-2$ vertices, and let $X:=\set{v_1,v_2}$ be the vertices outside $F$. We distinguish two cases according to the number of facets of $F$.
\begin{enumerate}   
\item \textbf{$F$ has at least $d+2$ facets}. We apply the induction hypothesis together with \cref{lem:facet-k-faces} to get the result. 
\item  \textbf{$F$ has exactly $d+1$ facets}. Here $f_{k}(F)\ge \tau_{k}(2(d-1)+2,d-1)$ by \cref{lem:lower-bound-dplus2-facets}.  Each vertex in $X$ is the apex of a pyramid over a $(d-2)$-face of $F$.  By analysing the $(d-2)$-faces of $F$ (see   \cref{rem:dplus2facets-facets}), we conclude that each vertex in $X$ has degree  at least $d+2$ in some facet of $P$ and at least $d+3$ in $P$.

We use \cref{cor:number-faces-outside-facet-practical} to find  faces $F_1,F_2$ in $P$ such that each $F_i$ has dimension $d-i+1$, contains $v_i\in X$, and excludes any $v_j\in X$ with $j<i$. Note that $F_1=P$. Hence the number of $k$-faces of $P$ that contain at least one vertex in $X$ is bounded below by 
\begin{equation*}
    \sum_{i=1}^2 f_{k-1}(F_{i}/v_{i})\ge \theta_{k-1}(d+3,d-1)+\theta_{k-1}(d+2,d-2).
\end{equation*}
We use \cref{cor:more-dpluss} to bound $f_{k-1}(F_{i}/v_{i})$. These elements settle the subcase. 
\end{enumerate}  
We now proceed with the details.

\begin{claim}
\label{cl:2d+2-dplus3-1} If $P$ has a facet $F$ with $2(d-1)+2=2d$ vertices, then $f_{k}(P)\ge \eta_{k}(2d+2,d)$. 
\end{claim}
\begin{claimproof} We consider two cases according to the number of $(d-2)$-faces in $F$.

\textbf{Suppose that $f_{d-2}(F)\ge d-1+3$}. For $k\in[2\ldots d-2]$, combining the induction hypothesis on $d-1$ and \cref{lem:facet-k-faces} gives that 
\begin{equation}
    \label{eq:Claim1_1}
    \begin{aligned}
f_{k}(P)&\ge f_{k}(F)+f_{k-1}(F)\\
&\ge\brac*{\binom{d}{k+1}+2\binom{d-1}{k+1}-\binom{d-3}{k+1}}+\brac*{\binom{d}{k}+2\binom{d-1}{k}-\binom{d-3}{k}}\\
&=\binom{d+1}{k+1}+2\binom{d}{k+1}-\binom{d-2}{k+1}= \eta_k(2d+2,d).
\end{aligned}
\end{equation}

\textbf{Assume that $f_{d-2}(F)= d-1+2$}. It follows that $F$ is  a $(d-1-a)$-fold pyramid over $T(m)\times T(a-m)$   for some $2\le a\le d-1$ and $2\le m\le \floor{a/2}$ (\cref{lem:dplus2facets}).  There are two vertices outside $F$, say $v_{1}$ and $v_{2}$. A facet $F_{2}$ containing $v_{2}$ but not $v_{1}$ must be a pyramid over a $(d-2)$-face $R$ of $F$; the same applies to a facet containing $v_{1}$ but not $v_{2}$. According to   \cref{rem:dplus2facets-facets}, for any $(d - 2)$-face $R$ of $F$, there are
just three possible values for $f_0(R)$, namely $2d - (a - m + 1)$ (which is $\ge d + 2$),
$2d-(m+1)$ (which is $\ge d+3$), and  $2d-1$.  This implies that $\deg_{F_{2}}(v_{2})\ge d+2$ and $\deg_{P}(v_{1}),\deg_{P}(v_{2})\ge d+3$. Consequently, thanks to  \cref{prop:number-faces-outside-facet}, the $k$-faces of $P$ include  the $k$-faces in $F$, at least $\tau_{k}(2(d-1)+2,d-1)$ by \cref{lem:lower-bound-dplus2-facets}, the $k$-faces in $P$ containing $v_{1}$, at least $\theta_{k-1}(d-1+4,d-1)$ by \cref{cor:more-dpluss}, plus the $k$-faces containing $v_{2}$ inside $F_{2}$, at least $\theta_{k-1}(d-2+4,d-2)$ by \cref{cor:more-dpluss}. When $k\ge 1$, this analysis yields the following: 
\begin{equation}
\label{eq_Claim1_2}
\begin{aligned}
f_{k}(P)&\ge f_{k}(F)+\sum_{i=1}^{2}\theta_{k-1}(d-i+4,d-i)\\
&\ge \tau_{k}(2(d-1)+2,d-1)+\brac*{\binom{d}{k}+\binom{d-1}{k}-\binom{d-4}{k}}\\
&\quad +\brac*{\binom{d-1}{k}+\binom{d-2}{k}-\binom{d-5}{k}}\\
&=\eta_{k}(2d+2,d)+\binom{d-2}{k+1}-\binom{\ceil{d/2}-1}{k+1}
-\binom{\ceil{d/2}-2}{k+1}-\binom{d-4}{k}\\
&\quad -\binom{d-5}{k}\\
&\ge \eta_{k}(2d+2,d).
\end{aligned}
\end{equation}

From \cref{lem:combinatorial-ineq}(iii), we have that $\binom{d-2}{k+1}-\binom{\ceil{d/2}-1}{k+1}
-\binom{\ceil{d/2}-2}{k+1}-\binom{d-4}{k}-\binom{d-5}{k}>0$ when $k\le d-3$ and that $\binom{d-2}{k+1}-\binom{\ceil{d/2}-1}{k+1}
-\binom{\ceil{d/2}-2}{k+1}-\binom{d-4}{k}-\binom{d-5}{k}=0$ for $k=d-2$. This settles this case.
\end{claimproof}

\subsection{Claim \ref{cl:2d+2-dplus3-2}: A facet with between $d$ and $2d-1$ vertices} We begin with an overview of the proof. Let $F$ be facet of $P$ with $f_{0}(F)=(d-1)+r$ for some $r\in [1\ldots d]$. If $1\le r\le d-1$, then we apply \cref{thm:at-most-2d} to obtain $f_{k}(F)\ge \theta_k(d-1+r,d-1)$. If $r=d$ and  $f_{d-2}(F)\ge d-1+3$, then we use $f_{k}(F)\ge \eta_{k}(2(d-1)+1,d-1)$ by \cref{thm:2dplus1}. If instead $r=d$ and  $f_{d-2}(F)= d+1$, we invoke \cref{lem:lower-bound-dplus2-facets} to get $f_{k}(F)\ge \tau_{k}(2(d-1)+1,d-1)$.

Let $X:=\set{v_{1},\ldots, v_{d-r+3}}$ be the set of vertices outside $F$. We apply \cref{cor:number-faces-outside-facet-practical} to find  faces $F_1,\ldots,F_{d-r+3}$ in $P$ such that each $F_i$ has dimension $d-i+1$, contains $v_i\in X$, and excludes any $v_j\in X$ with $j<i$. Thus, the number of $k$-faces of $P$ that contain at least one vertex in $X$ is bounded below by 
\begin{equation}
\label{eq:faces-in-X}
    \sum_{i=1}^{d-r+3} f_{k-1}(F_{i}/v_{i}).
\end{equation}

By \cref{lem:2d+s-dplus3-simple-pyramid}, we could assume that $P$ contains a vertex that is both nonpyramidal and nonsimple. Let  $v_1$ be  such a  vertex of maximum degree, and let $F$ be a facet not containing $v_1$.
The main difficulty with this approach lies in bounding $f_{k-1}(P/v_1)$, for which we only have the lower bound $\theta_{k-1}(d+1,d-1)$. The lower  bounds for $f_k(F)$ and $f_{k-1}(P/v_1)$, combined with  \eqref{eq:faces-in-X},  are insufficient to get the desired inequality. Therefore, we refine the argument by applying \cref{cor:number-faces-outside-facet-practical} twice: first to vertices inside a second facet $J_1\ne F$, and then to vertices outside $F\cup J_1$. For this approach to yield a larger lower bound  than \eqref{eq:faces-in-X}, we need at least two vertices outside $F\cup J_1$. We now proceed with the details.

\begin{claim} If $P$ has facets $F$ and $J_1$ such that $f_{0}(F)\le 2d-1$ and there are at least two vertices outside $F\cup J_1$, then $f_{k}(P)\ge \eta_{k}(2d+2,d)$.      
\label{cl:2d+2-dplus3-2}
\end{claim}
 \begin{claimproof} Assume that $f_{0}(F)=(d-1)+r$ for $r\in [1\ldots d]$ and let $X:=\set{v_{1},\ldots, v_{d-r+3}}$ be the set of vertices outside $F$. Additionally, let $J_{1}\cap X=\set{v'_{1},\ldots,v'_{t}}$.  Then $1\le t\le \card X-2=d-r+1$.  By \cref{cor:number-faces-outside-facet-practical}, there is a sequence $F'_1,\ldots, F'_{t}$ of faces within $J_{1}$ such that each $F'_i$ has dimension $d-1-i+1$ and contains $v'_i$ but does not contain any $v'_j$ with $j<i$.  The number of $k$-faces of $J_{1}$ that contain at least one of the vertices in $\set{v'_{1},\ldots,v'_{t}}$ is bounded from below by  	
\begin{align*}
\sum_{i=1}^{t}f_{k-1}(F'_{i}/v'_{i})&\ge \sum_{i=1}^{t} \binom{d-i}{k}.
\end{align*}
There are still $\card X-t=d-r+3-t\ge 2$ vertices outside $F\cup J_{1}$. Thus, by \cref{cor:number-faces-outside-facet-practical},  we have the following:

\begin{align}
\label{eq:2d+2-dplus3-2}
f_{k}(P)&\ge f_{k}(F)+\sum_{i=1}^{t}\binom{d-i}{k}+\sum_{j=1}^{d-r+3-t}\binom{d+1-j}{k}.
\end{align}

\textbf{In the case of $1\le r\le d-1$}, from \cref{thm:at-most-2d}, we get that 
\begin{equation}
\label{eq:atmost2d-dplus2-1}
f_{k}(F)\ge \binom{d}{k+1}+\binom{d-1}{k+1}-\binom{d-r}{k+1}.
\end{equation}
Combining \eqref{eq:2d+2-dplus3-2} and \eqref{eq:atmost2d-dplus2-1}, we get  
\begin{align*}
f_{k}(P)&\ge\brac*{\binom{d}{k+1}+\binom{d-1}{k+1}-\binom{d-r}{k+1}}+\sum_{i=1}^{t}\binom{d-i}{k}+\sum_{j=1}^{d-r+3-t}\binom{d+1-j}{k}.
\end{align*}
Furthermore, for $k\ge 2$, we can use  \cref{lem:combinatorial-identities} to write this expression as
\begin{equation}
\label{eq:2d+2-dplus3-3a}
\begin{aligned}
&=\eta_{k}(2d+2,d)-\sum_{\ell=1}^{d-r-1}\binom{d-r-\ell}{k}-\binom{d-2}{k}+\sum_{i=1}^{t-1}\binom{d-1-i}{k}\\
&\quad+\sum_{j=1}^{d-r+1-t}\binom{d-1-j}{k}.
\end{aligned}
\end{equation}

\textbf{In the case of $r=d$}, if $f_{d-2}(F)\ge d-1+3$, \cref{thm:2dplus1} gives 
\begin{equation}
\label{eq:2dplus1-dplus3}
f_{k}(F)\ge \eta_{k}(2(d-1)+1,d-1)=\binom{d}{k+1}+2\binom{d-1}{k+1}-\binom{d-2}{k+1}.
\end{equation}
If instead $f_{d-2}(F)=d-1+2$, then  \cref{lem:lower-bound-dplus2-facets} give that
\begin{equation}
\label{eq:2dplus1-dplus2}
\begin{aligned}
f_{k}(F)&\ge \tau_{k}(2(d-1)+1,d-1)\\
&=\binom{d}{k+1}+\binom{d-1}{k+1}+\binom{d-2}{k+1} -\binom{\ceil{(d-1)/2}}{k+1}-\binom{\ceil{(d-1)/2}-1}{k+1}.
\end{aligned}
\end{equation}

We analyse \eqref{eq:2d+2-dplus3-3a} for all $1\le r\le d-1$, and treat the case $r= d$   by combining \eqref{eq:2d+2-dplus3-2}  with either   \eqref{eq:2dplus1-dplus3} or \eqref{eq:2dplus1-dplus2}.

\begin{case} $1\le r\le d-1$ and $2\le t\le d-1$ (recall that $t\le d-r+1$). 
\label{case:2d+2-dplus3-1}
\end{case}
Here \eqref{eq:2d+2-dplus3-3a} becomes 
\begin{align}
\label{eq:2d+2-dplus3-4}
f_{k}(P)&\ge\eta_{k}(2d+2,d)+\brac*{\sum_{j=1}^{d-r+1-t}\underbrace{\binom{d-1-j}{k}-\binom{d-r-j}{k}}_{\ge 0\; \text{(for  $r\ge 1$)}}}\\
&\nonumber\quad+\brac*{\sum_{i=1}^{t-2}\underbrace{\binom{d-2-i}{k}-\binom{t-1-i}{k}}_{\ge 0\; \text{(as $t\le d-1$)}}}\\
&\nonumber\ge \eta_{k}(2d+2,d).
\end{align}
Equality holds precisely in the following scenarios: when $t=d-1$ and $r\in\{1,2\}$ (for all $k$); when $t=2$ and $r=1$ (for all $k$); and when $r=d+1-t$, $t\in [2\ldots d-1]$, and $k=d-2$.



\begin{case} $2\le r\le d-1$ and $t=1$. 
\label{case:2d+2-dplus3-2}
\end{case}
Here \eqref{eq:2d+2-dplus3-3a} becomes 

\begin{align}
\label{eq:2d+2-dplus3-5}
f_{k}(P)&\ge\eta_{k}(2d+2,d)+\brac*{\sum_{j=1}^{d-r-1}\underbrace{\binom{d-2-j}{k}-\binom{d-r-j}{k}}_{\ge 0\; \text{(for  $r\ge 2$)}}}\\
&\nonumber\ge \eta_{k}(2d+2,d).
\end{align}
Equality occurs only when $r=2$ (for all $k$), $r=d-1$ (for all $k$), or $k=d-2$ (for all $r$).


\begin{case} $r=d$.
\label{case:2d+2-dplus3-3a}
\end{case}

We must have that $t=1$. We consider two scenarios according to the number of $(d-2)$-faces in $F$. 

\textbf{First suppose that $f_{d-2}(F)\ge d-1+3$}. With the help of \cref{thm:2dplus1}, \eqref{eq:2d+2-dplus3-2} becomes
\begin{equation}
\label{eq:2d+2-claim2-case3-1}
\begin{aligned}
f_{k}(P)&\ge f_{k}(F)+\sum_{i=1}^{3}f_{k-1}(F_{i}/v_{i})\\
&\ge \brac*{\binom{d}{k+1}+2\binom{d-1}{k+1}-\binom{d-2}{k+1}}+\brac*{\binom{d}{k}+2\binom{d-1}{k}}\\
&=\binom{d+1}{k+1}+2\binom{d}{k+1}-\binom{d-2}{k+1}=\eta_{k}(2d+2,d).
\end{aligned}
\end{equation}


\textbf{Assume now that $f_{d-2}(F)= d-1+2$}. It follows that $F$ is  a $(d-1-a)$-fold pyramid over $T(m)\times T(a-m)$   for some $2\le a\le d-1$ and $2\le m\le \floor{a/2}$ (\cref{lem:dplus2facets}). According to   \cref{rem:dplus2facets-facets}, for any $(d-2)$-face $R$ contained in $F$, there are only three possible values for $f_0(R)$:  $2d-1-(a-m+1)$ (which is $\ge d+1$), $2d-1-(m+1)$ (which is $\ge d+2$), or $2d-2$.  Since $t=1$, $J_1$ must be a pyramid intersecting $F$ in a $(d-2)$-face, which has at least $d+1$ vertices. Thus $J_1$ has at least $d+2$ vertices; that is, $f_{0}(J_{1})=d-1+r'$ where $3\le r'\le d$.  

For $3\le r'\le d-1$, if we swap the roles of $J_{1}$ and $F$ and assume that $F$ contains $t'$ vertices from the vertices outside $J_{1}$, then $2\le t'\le d-2$. Consequently, the relevant pairs $(r',t')$ are covered in \cref{case:2d+2-dplus3-1}.  

It remains to consider $r'=d$. In this situation, both $F$ and $J_{1}$ are pyramids. The vertex in $J_{1}\setminus F$, say $v_{1}$, has degree $2d-2$ in $J_{1}$, and thus, degree at least $2d-1$ in $P$.  Let   $v_{2}$ and $v_{3}$ be the other two vertices outside $F$. By \cref{cor:number-faces-outside-facet-practical}, there is a sequence $F_1, F_2, F_{3}$ of faces of $P$ such that each $F_i$ has dimension $d-i+1$ and contains $v_i$, but does not contain any $v_j$ with $j<i$.

If the number of $(d-2)$-faces in $F_{1}/v_{1}$ is at least $(d-1)+3$ , then  \cref{thm:2dplus1}  yields, for $k\ge 2$, that
\begin{align*}
\sum_{i=1}^{3}f_{k-1}(F_{i}/v_{i})&\ge \eta_{k-1}(2(d-1)+1,d-1)+\binom{d-1}{k}+\binom{d-2}{k}.
\end{align*}
In this case, \cref{lem:lower-bound-dplus2-facets}  gives that 
\begin{equation}
\label{eq:2d+2-claim2-case3-2}
\begin{aligned}
f_{k}(P)&\ge f_{k}(F)+\sum_{i=1}^{3}f_{k-1}(F_{i}/v_{i})\\
&\ge \tau_{k}(2(d-1)+1,d-1)+\eta_{k-1}(2(d-1)+1,d-1)+\binom{d-1}{k}+\binom{d-2}{k},
\end{aligned}
\end{equation}
and \cref{lem:combinatorial-ineq}(iv) gives that 
\begin{equation*}
f_{k}(P)>\eta_{k}(2d+2,d).
\end{equation*}

Suppose the number of $(d-2)$-faces in $F_{1}/v_{1}$ is $(d-1)+2$. Then \cref{lem:lower-bound-dplus2-facets} yields, for $k\ge 2$, that
\begin{align*}
\sum_{i=1}^{3}f_{k-1}(F_{i}/v_{i})&\ge \tau_{k-1}(2(d-1)+1,d-1)+\binom{d-1}{k}+\binom{d-2}{k}.
\end{align*}
Thus  
\begin{equation}
\label{eq:2d+2-claim2-case3-3}
\begin{aligned}
f_{k}(P)&\ge f_{k}(F)+\sum_{i=1}^{3}f_{k-1}(F_{i}/v_{i})\\
&\ge \tau_{k}(2(d-1)+1,d-1)+\tau_{k-1}(2(d-1)+1,d-1)+\binom{d-1}{k}+\binom{d-2}{k}.
\end{aligned}
\end{equation}
\cref{lem:combinatorial-ineq}(v) now gives that 
\begin{equation*}
f_{k}(P)> \eta_{k}(2d+2,d).
\end{equation*}

\begin{case} $r=1$ and $t=d$. 
\label{case:2d+2-dplus3-3b}
\end{case}
In this case, $F$ is a simplex and $\card X=d+2$. If the facet $J_{1}$ is not a simplex, then $\dim(J_{1}\cap F)\ge 0$ and $d+1\le f_{0}(J_1)\le 2d-1$. Thus, by swapping the roles of $J_{1}$ and $F$, we have between 1 and $d-1$ vertices in $F\setminus J_{1}$. As a result,  the computation for $(r,t)$ where $2\le r\le d$ and $1\le t\le d-1$ would apply, namely Cases~\ref{case:2d+2-dplus3-1}--\ref{case:2d+2-dplus3-3a}.  Hence we may assume that $J_{1}$ is a simplex. It follows that $J_{1}$ is disjoint from $F$. 

Let $w_{1}$ and $w_{2}$ be the two vertices outside $F\cup J_{1}$.   Let $F_{1}$ be a facet containing $w_{1}$ but not $w_{2}$. Then $f_{0}(F_{1})\le 2d-1$. If $F_{1}$ is not a simplex, then   $J_{1}$ must share vertices with $F_{1}$, and thus $J_1\setminus F_1$ contains between $1$ and $d-1$ vertices. Additionally, there are (at least) two  vertices outside $F_{1}\cup J_{1}$. Thus, if $F_{1}$ played the role of $F$, then the computation for $(r,t)$ where $2\le r\le d$ and $1\le t\le d-1$ would apply. Hence we may assume that $F_{1}$ is a simplex. 

If $0\le \dim (F_{1}\cap F)< d-2$, then $J_{1}$ would share vertices with $F_{1}$, and $J_1\setminus F_1$ would contain between $2$ and $d-1$ vertices. Again,  if $F_{1}$ played the role of $F$, then the computation for $(r,t)$ where $r=1$ and $2\le t\le d-1$ would apply. As a consequence, we may assume that $F_{1}$ is a simplex intersecting $F$ in  either a ridge of $P$ or the empty set. By symmetry, we may also assume $F_{1}$ intersects $J_{1}$ in either a ridge of $P$ or the empty set. The same reasoning applies to the facets in $P$ containing $w_{2}$ but not $w_{1}$. Let $F_{2}$ be a facet containing $w_{2}$ but not $w_{1}$.  

Suppose that both $F_{1}$ and $F_{2}$ intersect $F$ in a ridge of $P$. Then the number of $k$-faces in $P$ is at least the sum of the following contributions: the number of $k$-faces of $F_{i}$ containing $w_{i}$, namely $\binom{d-1}{k}$; the number of $k$-faces of $F$,  namely $\binom{d}{k+1}$; and the number of $k$-faces in $P$ containing at least one vertex in $J_{1}$, which is at least $\sum_{i=1}^{d} \binom{d+1-i}{k}=\binom{d+1}{k+1}$ by \cref{prop:number-faces-outside-facet}. Thus 
\begin{equation}
\label{eq:2d+2-claim2-case4-1}
\begin{aligned}
f_{k}(P)&\ge 2\binom{d-1}{k} + \binom{d}{k+1}+\binom{d+1}{k+1}=\eta_{k}(2d+2,d)+\binom{d-1}{k}-\binom{d-2}{k}\\
&> \eta_{k}(2d+2,d).
\end{aligned}
\end{equation}
An analogous argument applies if both $F_{1}$ and $F_{2}$ intersect $J_1$ in a ridge. We may therefore assume that $F_1$ intersects $F$ in a ridge and $F_2$ intersects $J_1$ in a ridge. If  there were another facet $F_1'$ containing $w_1$ but not $w_2$ that intersects $J_1$ in a ridge, then both $F_2$ and $F_1'$ would intersect $J_1$ in a ridge. Hence, without loss of generality,  we may now assume that every facet containing $w_{1}$ but not $w_{2}$ intersects $F$ in a ridge, and every facet containing $w_{2}$ but not $w_{1}$ intersects $J_{1}$ in a ridge.

If there were no facet containing both $w_{1}$ and $w_{2}$, then there would be precisely $d$ facets containing $w_{1}$ and the collection of these facets and $F$  would form the boundary complex of a $d$-simplex, a contradiction. 

Let $F_{12}$ be a facet containing $w_{1}$ and $w_{2}$. If $f_0(F_{12})=2d$, then we find a scenario settled by \cref{cl:2d+2-dplus3-1}. Thus, we may assume that $f_{0}(F_{12})\le 2d-1$. Without loss of generality, we may further assume that $\dim(F_{12}\cap J_{1})\ge \dim(F_{12}\cap F)\ge  -1$, which implies that $\dim(F_{12}\cap J_{1})\ge 1$ for $d\ge 6$. When $\dim (F_{12}\cap F)=d-2$, we find that $f_{0}(F_{12})=2d$. Thus, we may further assume that $\dim (F_{12}\cap F)\le d-3$. Hence there are at least two vertices outside $F_{12}\cup J_{1}$, vertices in $F\setminus F_{12}$.

Suppose first that $F_{12}$ is not a simplex. If $F_{12}$ played the role of $F$, then $J_{1}\setminus F_{12}$ would contain between $1$ and $d-2$ vertices, and  the computation for $(r,t)$ where $ 2\le r\le d$ and $1\le t\le d-2$ would apply. Now suppose that  $F_{12}$ is a simplex. If $F_{12}$ played the role of $F$, then $J_{1}\setminus F_{12}$ would contain between $2$ and $d-2$ vertices, and the computation for $(r,t)$ where $r=1$ and $2\le t\le d-2$ would apply.  In either situation we are reduced to previously treated cases, which completes the proof of this case.

\begin{case} $r=1$ and $t=1$.
\label{case:2d+2-dplus3-5}
\end{case}

 In this case, $F$ is a simplex and $\card X=d+2$. Consider a facet $N_{1}$ with $d-1+r'$ vertices that  does not contain the vertex $v_{1}\in X$ but does contain some vertex in $X\setminus \set{v_{1}}$, say $v_{2}$. If $N_{1}$ contained $t'$ vertices from $X\setminus \set{v_{1}}$ with $2\le t'\le d$, then there would be at least two vertices outside $N_{1}\cup F$. If $N_{1}$ played the role of $J_{1}$, we would then be in one of the previously treated cases with  $(r,t)=(1,t')$. Therefore we may assume that $N_{1}$ contains either exactly one vertex from $X$ (namely $v_{2}$) or every vertex in $X\setminus \set{v_{1}}$.

\textit{Suppose that $N_{1}$ contains only $v_{2}$ from $X$.} In this case both $N_{1}$ and $F$ are simplices, and $N_{1}\cap F$ is a $(d-2)$-face. Thus we may take $N_1$ as the facet $J_1$. Now consider a facet $N_{2}$ different from $F$ that  does not contain $v_{2}\in X$. Then $N_2$ must  contain some vertex in $X\setminus \set{v_{2}}$, say $v_{p}$. By reasoning as before,  we may assume that $N_{2}$ contains either exactly one vertex from $X$ (namely $v_{p}$) or every vertex in $X\setminus \set{v_{2}}$. First suppose that $\dim (N_{2}\cap F)=d-2$. If $N_{2}$ contains only $v_p$ from $X$, then, by \cref{prop:number-faces-outside-facet}, the number of $k$-faces of $P$ containing at least one vertex in $X\setminus \set{v_{p},v_{2}}$ is at least $\sum_{i=1}^{d}\binom{d+1-i}{k}=\binom{d+1}{k+1}$ (since $k\ge 1$). The number of $k$-faces of $N_{1}$ containing $v_{2}$ is $\binom{d-1}{k}$ and the number of $k$-faces of $N_{2}$ containing $v_{p}$ is $\binom{d-1}{k}$. Thus 
\begin{equation}
\label{eq:2d+2-claim2-case5-1}
\begin{aligned}
f_{k}(P)&\ge f_{k}(F)+2\binom{d-1}{k}+\binom{d+1}{k+1}=\eta_{k}(2d+2,d)+\binom{d-1}{k}-\binom{d-2}{k}\\
&>\eta_{k}(2d+2,d).
\end{aligned}
\end{equation}
Moreover, if $N_{2}$ contains every vertex in $X\setminus \set{v_{2}}$, then	$f_{0}(N_{2})=d+1+d-1=2d$, which reduces to  \cref{cl:2d+2-dplus3-1}. Hence we may assume that such a facet $N_{2}$ does not intersect $F$ in a ridge of $P$, which implies that every facet intersecting $F$ in a ridge of $P$ contains $v_{2}$. However, this is impossible: dually, this would mean that every neighbour of the vertex $v_{F}$ conjugate to $F$ lies in the facet $K_{v_{2}}$ conjugate to $v_{2}$, but the facet $K_{v_{2}}$ is not the base of a pyramid in the dual polytope of $P$. This settles this scenario.

\textit{We may now assume that $N_{1}$ contains every vertex in $X\setminus \set{v_{1}}$.} If $N_{1}\cap F$ is a $(d-2)$-face, then $f_{0}(N_{1})=d+1+d-1=2d$, which reduces to  \cref{cl:2d+2-dplus3-1}. Thus we may
assume that $N_{1}\cap F$ is not a $(d-2)$-face. Consequently, we may finally assume that every facet intersecting $F$ in a ridge of $P$ contains $v_{1}$; however, this is impossible by the preceding argument. This contradiction settles \cref{cl:2d+2-dplus3-2}. 
\end{claimproof}
 
\subsection{Final part for the inequality statement of the theorem.} 
\label{subsec:final-part-inequality}

  By virtue of  \cref{cl:2d+2-dplus3-1}, we may assume that $P$ has no facet with $2d$ vertices. We claim that the intersection of any two facets of $P$ must be a ridge. This is easy to prove if one of them is the base of a pyramid, so let $F$ and $F_1$ be any two facets of $P$ such that $f_0(F)\le 2d-1$ and $f_0(F_1)\le 2d-1$.  By virtue of \cref{cl:2d+2-dplus3-2}, we may also assume that $F\cup F_1$ misses at most one vertex of $P$.  Suppose that $\dim (F_{1}\cap F)\in [-1\ldots d-3]$.  Then there is a facet $K$ other than $F_{1}$ and $F$ that contains $F_{1}\cap F$. Furthermore,  we see that $K$ cannot miss two vertices in $F\setminus F_1$ and similarly cannot miss two vertices  in $F_1\setminus F$.  This implies that there is a unique vertex in $F\setminus K$ and a unique vertex in $F_1\setminus K$, and no vertices at all outside $F\cup F_1\cup K$. Thus $K$ has $2d$ vertices in total, contrary to our assumption. Hence $\dim (F_{1}\cap F)= d-2$.
 
 Thus, $P$ is a  super-Kirkman  polytope, meaning that every pair of facets intersects in a ridge of the polytope \cite{PinWanYos24}.
Since $P$ is super-Kirkman and has at least $d+3$ facets, every facet of $P$ must contain at least $d+2$ ridges of $P$.  The smallest simple $(d-1)$-polytope with at least $d+2$ $(d-2)$-faces contains at least  $3d-4$ vertices  (\cref{thm:simple-lbt}), which is larger than $2d+1$ for $d\ge 6$. The following is now plain.

\begin{remark}
\label{rmk:facet-nonsimple-dplus3}
Every facet of $P$ is a nonsimple polytope: it has either between $d+1$ and $2d-1$ vertices or $2d+1$ vertices.
\end{remark}

If every vertex in $P$ is either pyramidal or simple, then $P$ is a multifold pyramid over a simple polytope. In this instance, the result follows from \cref{lem:2d+s-dplus3-simple-pyramid}. Thus we may assume that $P$ have a nonpyramidal, nonsimple vertex $v_{1}$. The vertex $v_{1}$ has degree at least $d+4$ in $P$, as every nonsimple vertex in a super-Kirkman $d$-polytope has degree at least $d+4$ \cite[Thm.~3.9]{PinWanYos24}.
Let $F$  be a facet in $P$ that does not contain $v_{1}$, and let $X$ be the set of vertices outside $F$.  From \cref{rmk:facet-nonsimple-dplus3}, we get that $f_{0}(F)=d-1+r$ with $r\in [2\ldots d]$. 

\textbf{Suppose that $r=d$}. There are three vertices outside $F$, say $v_{1}$, $v_{2}$, and $v_{3}$. By \cref{cor:number-faces-outside-facet-practical}, there is a sequence $F_1, F_2, F_{3}$ of faces of $P$ such that each $F_i$ has dimension $d-i+1$ and contains $v_i$, but does not contain any $v_j$ with $j<i$. The number of $k$-faces of $P$ that contain $v_{1}$ is at least $\theta_{k-1}(d+4,d-1)$ (by \cref{cor:more-dpluss}), so the number of $k$-faces that contain at least one vertex in $\set{v_{1},v_{2},v_{3}}$ is bounded below by  	
\begin{align*}
\sum_{i=1}^{3}f_{k-1}(F_{i}/v_{i})&\ge \theta_{k-1}(d+4,d-1)+\theta_{k-1}(d-1,d-2)+\theta_{k-1}(d-2,d-3)\\
&= \binom{d}{k}+\binom{d-1}{k}-\binom{d-5}{k}+\binom{d-1}{k}+\binom{d-2}{k}\\
&= \binom{d}{k}+2\binom{d-1}{k}+\sum_{i=1}^{3}\binom{d-2-i}{k-1}.
\end{align*}

In this case, \cref{thm:2dplus1} yields that 
\begin{equation}\label{eq:sk_r=d}
    \begin{aligned}
     f_{k}(P)&\ge f_{k}(F)+\sum_{i=1}^{3}f_{k-1}(F_{i}/v_{i})\\
     &\ge \brac*{\binom{d}{k+1}+2\binom{d-1}{k+1}-\binom{d-2}{k+1}}+\binom{d}{k}+2\binom{d-1}{k}\\
     &\quad +\sum_{i=1}^{3}\binom{d-2-i}{k-1}\\
     &=\binom{d+1}{k+1}+2\binom{d}{k+1}-\binom{d-2}{k+1}+\sum_{i=1}^{3}\binom{d-2-i}{k-1}\\
     &=\eta_{k}(2d+2,d)+\sum_{i=1}^{3}\binom{d-2-i}{k-1}\\
&>\eta_{k}(2d+2,d).
    \end{aligned}
\end{equation}

\textbf{Suppose that $3\le r\le d-1$}. There are $d-r+3$ vertices outside $F$, including $v_{1}$. By \cref{cor:number-faces-outside-facet-practical}, the number of $k$-faces of $P$ that contain at least one of the vertices outside $F$ is bounded from below by 
\begin{equation*}
\binom{d}{k}+\binom{d-1}{k}-\binom{d-5}{k}+\sum_{i=1}^{d-r+2}\binom{d-i}{k}.
\end{equation*}
Additionally, the facet $F$ is not a triplex, as it contains $d+2$ ridges of $P$ (\cref{thm:at-most-2d}). Thus $f_{k}(F)>\theta_{k}(d-1+r,d-1)$. We now have that
\begin{equation}\label{eq:sk_3<=r<=d-1}
\begin{aligned}
f_{k}(P)&\ge f_{k}(F)+\binom{d}{k}+\binom{d-1}{k}-\binom{d-5}{k}+\sum_{i=1}^{d-r+2}\binom{d-i}{k}\\
& > \brac*{\binom{d}{k+1}+\binom{d-1}{k+1}-\binom{d-r}{k+1}}+\binom{d}{k}+\binom{d-1}{k}-\binom{d-5}{k} \\
&\quad +\sum_{i=1}^{d-r+2}\binom{d-i}{k}\\
&=\eta_{k}(2d+2,d)-\binom{d-r}{k+1}-\binom{d-5}{k}+\sum_{i=3}^{d-r+2}\binom{d-i}{k}\\
&=\eta_{k}(2d+2,d)-\sum_{\ell=1}^{d-r-1}\binom{d-r-\ell}{k}-\binom{d-5}{k}+\sum_{i=3}^{d-r+2}\binom{d-i}{k}\\
&=\eta_{k}(2d+2,d)+\underbrace{\brac*{\sum_{\ell=1}^{d-r-1}\binom{d-3-\ell}{k}-\binom{d-r-\ell}{k}}}_{\ge 0\;\text{(for $r\ge 3$)}}\\
&\quad +\underbrace{\brac*{\binom{d-3}{k}-\binom{d-5}{k}}}_{\ge 0}\\
&\ge \eta_{k}(2d+2,d).
\end{aligned}
\end{equation}
Hence $f_{k}(P)>\eta_{k}(2d+2,d)$ for all $k\in [2\ldots d-2]$.

Finally, \textbf{assume that $r=2$}. A facet $F_{1}$ other than $F$ not containing $v_{1}$ must contain $\card X-1$ vertices from $X$ by \cref{cl:2d+2-dplus3-2}. Furthermore, $\dim (F_{1}\cap F)=d-2$, implying that $f_{0}(F_{1}\cap F)=d-1$ or $d$. In the former case,  $f_{0}(F_{1})=2d-1$, while in the latter case $f_{0}(F_{1})=2d$. Both scenarios have already been covered.  This completes the proof of the theorem.

\subsection{The equality statement of the theorem}

The equality part is  addressed  in \cref{cl:2d+2-dplus3-equality}, where we examine all instances in \cref{cl:2d+2-dplus3-1}, \cref{cl:2d+2-dplus3-2}, and \cref{subsec:final-part-inequality} in which the inequalities can hold with equality, relying on the analysis developed in the inequality part.

\begin{claim}
    Let $P$ be a $d$-polytope with $2d+2$ vertices, at least $d+3$ facets, and $f_k(P)=\eta_k(2d+2,d)$ for some $k\in [1\ldots d-2]$. Then $P$ has exactly $d+3$ facets. 
    \label{cl:2d+2-dplus3-equality}
\end{claim}

    \begin{claimproof}
We proceed by induction  on $d\ge 3$ for all $k\in [1\ldots d-2]$. The case $k=1$ was settled in \cite[Thm.~22]{PinUgoYos22} for all $d$. For $d=3, 4,5$, \cref{thm:2d+2-small} ensures that  equality holds for some $k\in [2\ldots d-2]$ if and only if $P$ has exactly $d+3$ facets.
      
Let $d\ge 6$, and let $P$ be a $d$-polytope with $2d+2$ vertices and at least $d+3$ facets such that  $f_k(P)=\eta_k(2d+2, d)$ for some $k\in [2\ldots d-2]$.  We examine the possible equality cases from the proof of \cref{thm:2dplus2-restated}, beginning with those in  \cref{cl:2d+2-dplus3-1,cl:2d+2-dplus3-2}.

\subsection*{(i) \cref{cl:2d+2-dplus3-1}:  $P$ has a facet with $f_{0}(P)-2$ vertices.} Let $v_1, v_2$ be the vertices outside $F$, and  consider two cases based on the number of $(d-2)$-faces in $F$.
      
\textbf{Suppose $F$ has at least $d+2$ $(d-2)$-faces}.  If $f_k(P)=\eta_k(2d+2, d)$ for some $k\in [2\ldots d-2]$ in \eqref{eq:Claim1_1}, then $f_k(F)=\eta_k(2d, d-1)$, $f_{k-1}(F)=\eta_{k-1}(2d, d-1)$, and each $k$-face of $P$ is either a $k$-face of $F$ or a $k$-face intersecting $F$ in a $(k-1)$-face. By induction,  $F$ has exactly $d+2$ $(d-2)$-faces. If there were a  facet $F'$ not intersecting $F$ in a ridge of $P$, then $F'$ would be a two-fold pyramid over a $(d-3)$-face of $F$. However, this would imply the existence of  a $k$-face in $P$ that is a two-fold pyramid over a $(k-2)$-face of $F$, contradicting our equality condition. Thus every facet of $P$ intersects $F$ in a $(d-2)$-face, which  implies that $P$ has exactly $d+3$ facets.

\textbf{Suppose $F$ has exactly $d+1$ $(d-2)$-faces}. If $f_k(P)=\eta_k(2d+2, d)$ for some $k\in [2\ldots d-2]$  in \eqref{eq_Claim1_2}, then $f_k(F)=\tau_k(2d,d-1)$;  the number of $k$-faces containing $v_1$ is exactly $\theta_{k-1}(d+3, d-1)$; and  there is a facet $F_2$ containing $v_2$ but not $v_1$, where $v_2$  lies in exactly $\theta_{k-1}(d+2, d-2)$ $k$-faces. We need the following basic fact, inspired by \cite[Sec.~4]{Xue21}; we omit its proof.

\begin{fact}
\label{fact:two-faces}
Let $F'$ and $F''$ be two distinct facets of a $d$-polytope that both contain a vertex $v$, and let $1\le k\le d-1$. Then there exists a $k$-face of $F'$ containing $v$ that is not a $k$-face of $F''$.
\end{fact}

If there were another facet $F_2'\ne F_2$ containing $v_2$ but not $v_1$, then, by \cref{fact:two-faces}, the number of $k$-faces containing $v_2$ but not $v_1$ would be strictly larger than $\theta_{k-1}(d+2, d-2)$, contradicting our equality condition. Thus,  apart from the facet $F$, $P$ has exactly $d+1$ facets containing $v_1$ (by \cref{cor:more-dpluss}) and a unique facet containing $v_2$ but not $v_1$. Hence $P$ has exactly $d+3$ facets.

\subsection*{(ii) \cref{cl:2d+2-dplus3-2}:   $P$ has facets $F$ and $J_1$ such that $f_{0}(F)\le 2d-1$ and there are at least two vertices outside $F\cup J_1$.} 

As in the proof of the claim, let $1\le r\le d$ and $1\le t\le d-r+1$, and assume  that $f_0(F)=(d-1)+r$. Additionally, let $X=\{v_1, v_2,\dots, v_{d-r+3}\}$ be the set of  vertices outside $F$, let $J_1\cap X=\set{v_1',\dots, v_t'}$, and let $X\setminus J_1=\set{v_1'',\dots, v_{d-r+3-t}''}$. We begin with a simple fact.

\begin{fact}
\label{fact:minimum-var}
   If the $d-r+3-t$ vertices in $X\setminus J_1$ are contained in precisely $\sum_{i=1}^{d-r+3-t}\binom{d+1-i}{k}$ $k$-faces of $P$, then the number of $k$-faces in $P$ that contain the vertex $v_i''$ and no vertex $v_j''$ with $j<i$ is precisely $\binom{d+1-i}{k}$.  
\end{fact}
\begin{factproof}  By \cref{cor:number-faces-outside-facet-practical}, there is a sequence $Y_1,\ldots, Y_{d-r+3-t}$ of faces in $P$ such that each $Y_i$ has dimension $d+1-i$ and contains $v''_i$ but does not contain any $v''_j$ with $j<i$. The minimum number of $k$-faces within the $(d+1-i)$-face $Y_i$ that contain the vertex $v''_{i}$ is attained when the vertex figure $Y_i/v''_i$ is  a $(d-i)$-simplex.  Since $Y_i$ excludes each $v''_j$ with $j<i$ and the number of $k$-faces of $P$ containing some vertex in $X\setminus J_1$ is precisely $\sum_{i=1}^{d-r+3-t}\binom{d+1-i}{k}$, each vertex figure must be a minimiser. See also \cref{fact:two-faces}.
\end{factproof}


Following the proof of \cref{cl:2d+2-dplus3-2}, there are five cases to consider based on the pair $(r,t)$. 

\begin{case}   $1\le r\le d-1$  and $2\le t\le d-1$.
\label{case:1-equality}
\end{case}
If $f_k(P)=\eta_k(2d+2,d)$ in \eqref{eq:2d+2-dplus3-2} (see also \eqref{eq:2d+2-dplus3-4}), then (a) the number of $k$-faces of $P$ containing some vertices in $J_1\cap X$ but no vertex in $X\setminus J_1$ is exactly $\sum_{i=1}^{t}\binom{d-i}{k}$; (b) every vertex in $J_1\cap X$ is simple in  $J_1$; (c) the number of $k$-faces of $P$ containing some vertex in $X\setminus J_1$ is exactly $\sum_{j=1}^{d-r+3-t}\binom{d+1-j}{k}$; and (d) every vertex in $X\setminus J_1$ is simple in $P$.

Since $v_1''$ is simple in $P$, it is contained in precisely $d$ facets of $P$. There is a unique facet  $I_2$ of $P$ containing  $v_2''$ but not $v_1''$. Indeed, suppose that there are two facets $I_2$ and $I_2'$ of $P$ containing  $v_2''$ but not $v_1''$. Then, by \cref{fact:two-faces},  $v_2''$  would contribute more than $\binom{d-1}{k}$ to $\sum_{j=1}^{d-r+3-t}\binom{d+1-j}{k}$, contradicting \cref{fact:minimum-var}.

If there were a facet of $P$  different from $ J_1$ that contains a vertex in $J_1\cap X$ but no vertex from $X\setminus J_1$, say $J_1'$, then we could select such a vertex as $v_1'$. The contribution of the $k$-faces containing $v_1'$ in $J_1\cup J_1'$  would then exceed $\binom{d-1}{k}$ (see \cref{fact:two-faces}), yielding more than $\sum_{i=1}^{t}\binom{d-i}{k}$ $k$-faces of $P$ containing a vertex in $J_1\cap X$ but no vertex in $X\setminus J_1$. Thus $P$ has exactly $d+3$ facets: the $d$ facets containing $v_1''$, together with $I_2$, $J_1$, and $F$.

\begin{case}
    $2\le r\le d-1$ and $t=1$.
\end{case}
If $f_k(P)=\eta_k(2d+2,d)$ in \eqref{eq:2d+2-dplus3-2} (see also \eqref{eq:2d+2-dplus3-5}), then the same conditions as in \cref{case:1-equality} apply, and therefore so does the same argument. Thus $P$ has exactly $d+3$ facets.

\begin{case}
   $r=d$ (and $1\le t\le d-r+1$).
\end{case}
We consider two subcases according to the number of $(d-2)$-faces of $F$. If $f_{d-2}(F)\ge d-1+3$  and  $f_k(P)=\eta_k(2d+2,d)$ in \eqref{eq:2d+2-dplus3-2} (see also \eqref{eq:2d+2-dplus3-5}), then the same conditions and argument  in \cref{case:1-equality} apply, and therefore $P$ has exactly $d+3$ facets.
 
Suppose $f_{d-2}(F)= d-1+2$. Then $F$ is  a $(d-1-a)$-fold pyramid over $T(m)\times T(a-m)$   for some $2\le a\le d-1$ and $2\le m\le \floor{a/2}$ (\cref{lem:dplus2facets}). Since $t=1$, $J_1$ must be a pyramid intersecting $F$ in a $(d-2)$-face, which has at least $d+1$ vertices because $m\ge 2$ (\cref{rem:dplus2facets-facets}). Thus $f_{0}(J_{1})=d-1+r'$ where $3\le r'\le d$.  For $3\le r'\le d-1$, if we swap the roles of $J_{1}$ and $F$ and assume that $F$ contains $t'$ vertices from the vertices outside $J_{1}$, then $2\le t'\le d-2$. Consequently, the pairs $(r',t')$ are covered in \cref{case:2d+2-dplus3-1}. If $r'=d$, then we have strict inequalities in all the scenarios; see \eqref{eq:2d+2-claim2-case3-2}  and \eqref{eq:2d+2-claim2-case3-3}.

\begin{case}
   $r=1$ and $t=d$.
\end{case}
Here $F$ is a simplex. If the facet $J_{1}$ were not a simplex, then $\dim(J_{1}\cap F)\ge 0$ and $d+1\le f_{0}(J_1)\le 2d-1$. Thus, by swapping the roles of $J_{1}$ and $F$, there would have between 1 and $d-1$ vertices in $F\setminus J_{1}$, reducing this scenario to Cases~\ref{case:2d+2-dplus3-1}--\ref{case:2d+2-dplus3-3a}. If $J_1$ is a simplex, then we are reduced to \cref{cl:2d+2-dplus3-1}, to a previously treated case in \cref{cl:2d+2-dplus3-2}, or we obtain  strict inequality in \eqref{eq:2d+2-claim2-case4-1}.

\begin{case}
   $r=1$ and $t=1$.
\end{case}
Here $F$ is also a simplex and $J_1$ is a pyramid. In this case, the analysis yields strict inequality \eqref{eq:2d+2-claim2-case5-1}, a reduction to \cref{cl:2d+2-dplus3-1}, or a reduction to a previously treated case in \cref{cl:2d+2-dplus3-2}. 

This concludes the analysis of the equality cases arising from \cref{cl:2d+2-dplus3-1,cl:2d+2-dplus3-2}. We now deal with cases  considered in \cref{subsec:final-part-inequality}. 

By \cref{cl:2d+2-dplus3-1}, we may assume that $P$ has no facet with $2d$ vertices. We also know that $P$ is a super-Kirkman polytope. Additionally, every facet of $P$ is nonsimple; and its number of vertices is either $2d+1$ or in the range $[d+1\ldots 2d-1]$ (see \cref{rmk:facet-nonsimple-dplus3}). 
      
If every vertex in $P$ is either pyramidal or simple, then $P$ is a multifold pyramid over a simple polytope. By \cref{lem:2d+s-dplus3-simple-pyramid}, $P$ has exactly $d+3$ facets. Suppose now that there is a vertex $v_1$ that is nonpyramidal and nonsimple, and let $F$ be a facet that does not contain $v_1$, where $f_0(F)=d-1+r$ with $r\in [2\ldots d]$.
If $r=d$,  we have strict inequality in \eqref{eq:sk_r=d}, while if $3\le r\le d-1$, we have strict inequality in \eqref{eq:sk_3<=r<=d-1}. If $r=2$, then there is a facet $F_{1}\ne F$ not containing $v_{1}$ with $f_{0}(F_{1})\in \set{2d-1,2d}$. Both scenarios have already been covered. This completes the proof of the claim.
\end{claimproof}

The proof of \cref{thm:2dplus2-restated} is now complete.

\section{Combinatorial identities and proofs}
We start with a combinatorial lemma whose three first expressions are from \cite{Xue21}.
\begin{lemma}[Combinatorial equalities and inequalities] For all integers $d\ge 2$, $k\in [1\ldots d-1]$, the following hold
\label{lem:combinatorial-identities}
\begin{enumerate}[{\rm (i)}] 
\item If $2\le r\le s\le d$ are integers, then $\theta_{k}(d+s-r,d-1)+\sum_{i=1}^{r}\binom{d+1-i}{k}\ge \theta_{k}(d+s,d)$, with equality only if $r=2$ or $r=s$.
\item If $n\ge c$ are positive integers, then $\binom{n}{c}=\binom{n-1}{c-1}+\binom{n-1}{c}$.
\item If $n\ge c$ and $n\ge a$  are positive integers, then  $\binom{n}{c}-\binom{n-a}{c}=\sum_{i=1}^{a}\binom{n-i}{c-1}$.
\item If $n,c$ are positive integers, then  $\binom{n}{c}=\sum_{i=1}^{n}\binom{n-i}{c-1}$.
 \item (Vandermonde's identity) If $n,a,c$ are nonnegative integers, then \[\binom{n+a}{c}=
\sum_{i=0}^{c}\binom{n}{i}\binom{a}{c-i}.\]   
\end{enumerate}  
\end{lemma}
\begin{proof}

The proof of (i) is embedded in the proof of \cite[Thm.~3.2]{Xue21}; it is also spelled out in \cite[Claim 1]{Pin24}. Repeated applications of (ii) yield (iii) and (iv), while (v) is well known. \end{proof}

\begin{lemma} The following inequalities hold.
\begin{enumerate}[{\rm (i) }]
\item If  $d\ge 9$ and  $1\le k \le \ceil{d/3}-2$, then $\ff_{k}(2d+2,d)< \hh_{k}(2d+2,d)$.
\item If  $d\ge 9$ and  $\floor{0.4d}\le k \le d-1$, then $\ff_{k}(2d+2,d)> \hh_{k}(2d+2,d)$.
\item If  $d\ge 6$ and  $1\le k \le d-3$, then 
\begin{equation*}
\binom{d-2}{k+1}-\binom{\ceil{d/2}-1}{k+1}-\binom{\ceil{d/2}-2}{k+1}-\binom{d-4}{k}-\binom{d-5}{k}>0.
\end{equation*}
If $k=d-2$ or $d-1$, then the expression is zero.
\item If  $d\ge 4$ and  $2\le k \le d-2$, then 
\begin{equation*}
\ff_{k}(2d+2,d)< \hh_{k}(2(d-1)+1,d-1)+\ff_{k-1}(2(d-1)+1,d-1)+\binom{d-1}{k}+\binom{d-2}{k}.
\end{equation*}
\item If  $d\ge 5$ and  $1\le k \le d-2$, then 
\begin{equation*}
\ff_{k}(2d+2,d)< \hh_{k}(2(d-1)+1,d-1)+\hh_{k-1}(2(d-1)+1,d-1)+\binom{d-1}{k}+\binom{d-2}{k}.
\end{equation*}
\end{enumerate}
\label{lem:combinatorial-ineq}
\end{lemma}
\begin{proof} We reason as in \cite[Prop.~6.1]{Xue24}. 

(i) Write $\f_{k}(d):=\hh_{k}(2d+2,d)-\ff_{k}(2d+2,d)$. We can verify that $\f_{1}(9)=4$. 
 We prove the statement by induction on $d\ge 9$ for all $1\le k \le \ceil{d/3}-2$. Assume that $d\ge 10$. We have that 
\begin{align*}
\hh_{k}(2d+2,d)&=\binom{d + 1}{k + 1}+\binom{d}{k + 1}+\binom{d - 1}{k + 1}-\binom{\ceil{(d+1)/2}-1}{k + 1}\\
&\quad-  \binom{\ceil{(d+1)/2}-2}{k + 1}\\
\ff_{k}(2d+2,d)&=\binom{d + 1}{k + 1}  +  2\binom{d}{k + 1}-\binom{d - 2}{k + 1},\\
\intertext{which, after some simplifications, gives that}
\f_{k}(d)&=\binom{d -2}{k + 1}-\binom{d -1}{k}-\binom{\ceil{(d+1)/2}-1}{k + 1}-\binom{\ceil{(d+1)/2}-2}{k + 1}.
\end{align*}
Since 
\begin{align*}
\f_{k}(d-1)&=\binom{d -3}{k + 1}-\binom{d -2}{k}-\binom{\ceil{d/2}-1}{k + 1}-\binom{\ceil{d/2}-2}{k + 1}\\
\f_{k-1}(d-1)&=\binom{d -3}{k}-\binom{d -2}{k-1}-\binom{\ceil{d/2}-1}{k}-\binom{\ceil{d/2}-2}{k},
\end{align*}
from \cref{lem:combinatorial-identities}(ii) it follows that
\begin{equation}
\label{lem:combinatorial-ineq-1}
\f_{k}(d-1)+\f_{k-1}(d-1)=\binom{d -2}{k + 1}-\binom{d -1}{k}-\binom{\ceil{d/2}}{k + 1}-\binom{\ceil{d/2}-1}{k + 1}.
\end{equation}
If $d$ is even, then  $\f_{k}(d)=\f_{k}(d-1)+\f_{k-1}(d-1)$.	If $d=2p+1$ for some integer $p$, then
\begin{align}
\label{lem:combinatorial-ineq-2}
\f_{k}(d)-(\f_{k}(d-1)+\f_{k-1}(d-1))&= \binom{p-1}{k}+\binom{p}{k}> 0.
\end{align}
If both $k-1, k\le \ceil{(d-1)/3}-2$, combining \eqref{lem:combinatorial-ineq-1} and \eqref{lem:combinatorial-ineq-2} and the induction hypothesis on $d-1$ give that $\f_{k}(d)> 0$ . Thus, it remains to consider the case
\begin{equation*}
\ceil{(d-1)/3}-1\le k\le \ceil{d/3}-2,
\end{equation*}
which reduces to the pairs $(d,k)=(6p+1,2p-1), (6p+4,2p)$ for some integer $p\ge 1$. Consider the pair $(6p+1,2p-1)$.
\begin{align*}
\label{lem:combinatorial-ineq-3}
\f_{2p-1}(6p+1)&=\binom{6p-1}{2p}-\binom{6p}{2p-1}-\binom{3p}{2p}-\binom{3p-1}{2p}\\
&=\frac{p+1}{3p}\binom{6p}{2p-1}-4\binom{3p-1}{2p}\\
&=\frac{(3p-1)!}{(2p)!(p-1)!} \brac*{2\cdot\frac{3p+1}{p+2}\cdots \frac{6p}{4p+1}-4}\\
&=\binom{3p-1}{2p}\brac*{2\cdot\frac{3p+1}{p+2}\cdots \frac{6p}{4p+1}-4}.
\end{align*}	
From the last expression, we can gather that $\f_{2p-1}(6p+1)> 0$. 

Now consider the pair $(6p+4,2p)$.
\begin{align*}
\f_{2p}(6p+4)&=\binom{6p+2}{2p+1}-\binom{6p+3}{2p}-\binom{3p+2}{2p+1}-\binom{3p+1}{2p+1}\\
&=\frac{2p+3}{2(4p+3)}\binom{6p+2}{2p+1}-\frac{4p+3}{p+1}\binom{3p+1}{2p+1}\\
&=\binom{3p+1}{2p+1}\frac{4p+3}{p+1}\brac*{\frac{2p+3}{2(4p+3)}\cdot \frac{3p+2}{p+2}\cdots\frac{4p+2}{2p+2}\cdot\frac{4p+4}{2p+3} \cdots  \frac{6p+2}{4p+1}-1}\\
&=\binom{3p+1}{2p+1}\frac{4p+3}{p+1}\brac*{\frac{2p+3}{2p+3}\cdot \frac{1}{2}\cdot\frac{3p+2}{p+2}\cdots\frac{4p+2}{2p+2}\cdot\frac{4p+4}{4p+3} \cdots  \frac{6p+2}{4p+1}-1}.
\end{align*}
	
From the last expression, we can gather that $\f_{2p}(6p+4)> 0$. 

(ii)	Write $\f_{k}(d):=\ff_{k}(2d+2,d)-\hh_{k}(2d+2,d)$. 	After some simplification, we have \[\f_k(d)=\binom{d-1}{k}-\binom{d-2}{k+1}+\binom{\ceil{(d+1)/2}-1}{k+1}+\binom{\ceil{(d+1)/2}-2}{k+1}.\]
	
	To show that $\f_k(d)>0$ for $d\ge 9$, it suffices to prove that $\f'_k(d):=\binom{d-1}{k}-\binom{d-2}{k+1}>0$. For $k=d-2,d-1$, it is straightforward that $\f'_k(d)>0$. Thus, assume that $k\le d-3$. 
		\begin{align*}
			\f'_k(d)=&\binom{d-1}{k}-\binom{d-2}{k+1}\\
			=&\frac{(d-1)!}{(d-k-1)!k!}-\frac{(d-2)!}{(d-k-3)!(k+1)!}\\
			=&\frac{(d-2)!}{(d-k-3)!k!}\left(\frac{d-1}{(d-k-1)(d-k-2)}-\frac{1}{k+1}\right)\\
			=&\frac{(d-2)!}{(d-k-3)!k!}\left(\frac{(d-1)(k+1)-(d-k-1)(d-k-2)}{(d-k-1)(d-k-2)(k+1)}\right)\\
			=&\frac{(d-2)!}{(d-k-3)!k!}\left(\frac{-k^2+(3d-4)k+(-d^2+4d-3)}{(d-k-1)(d-k-2)(k+1)}\right).
		\end{align*}
		The expression $-k^2+(3d-4)k+(-d^2+4d-3)$ is positive when $k$ is between $k_1:=(3d-4-\sqrt{5d^2-8d+4})/2$ and $k_2:=(3d-4+\sqrt{5d^2-8d+4})/2$. Since $k_1<\floor{0.4d}$ and $k_2>d-1$ for all $d\ge 9$, it follows that $\f'_k(d)>0$ for $\floor{0.4d}\leq k\leq d-3$.

(iii) 	Write $\f_k(d):=\binom{d-2}{k+1}-\binom{\ceil{d/2}-1}{k+1}-\binom{\ceil{d/2}-2}{k+1}-\binom{d-4}{k}-\binom{d-5}{k}$. From \cref{lem:combinatorial-identities} we get that
\begin{equation}
\label{lem:combinatorial-ineq-4}
			\f_k(d)=\brac*{\binom{d-5}{k+1}-\binom{\ceil{d/2}-2}{k+1}}+\brac*{\binom{d-3}{k}-\binom{\ceil{d/2}-1}{k+1}}.
		\end{equation}
	
Suppose that $d=2p+1$ for some $p\geq 3$. Then 
\begin{align*}
			\f_k(d)=&\brac*{\binom{2p-4}{k+1}-\binom{p-1}{k+1}}+\brac*{\binom{2p-2}{k}-\binom{p}{k+1}}.
\end{align*}
For $k+1\geq p$, we have $\f_k(d)>0$. Hence assume that $k+1<p$. If $p\le 4$, then we can easily verify that $\f_k(d)>0$ for all $k$. So assume $p\ge 5$. On one hand, 
 \begin{align*}
			\binom{2p-4}{k+1}-\binom{p-1}{k+1}=\sum_{\ell=1}^{p-3}\binom{2p-4-\ell}{k} > (p-3)\binom{p-1}{k},
		\end{align*}        
and on the other hand, 
       \begin{align*} 
            \brac*{\binom{2p-2}{k}-\binom{p}{k+1}} > \binom{p-1}{k}-\binom{p}{k+1}=\binom{p-1}{k}-\frac{p}{k+1}\binom{p-1}{k}.
		\end{align*}
Thus, we have 
\begin{align*}
    \f_k(d)& > (p-2)\binom{p-1}{k}-\frac{p}{k+1}\binom{p-1}{k} \ge (p-2)\binom{p-1}{k}-\frac{p}{2}\binom{p-1}{k}\\
    &\ge \frac{p-4}{2}\binom{p-1}{k}>0.
\end{align*}
	
		Suppose that $d=2p$ for some $p\geq 3$. Then
\begin{align*}
\f_k(d)=&\brac*{\binom{2p-5}{k+1}-\binom{p-2}{k+1}}+\brac*{\binom{2p-3}{k}-\binom{p-1}{k+1}}.     \end{align*} 
For $k+2\geq p$, we have $\f_k(d)>0$. Hence assume that $k+2<p$. If $p\le 4$, then we can easily verify that $\f_k(d)>0$ for all $k$. So assume $p\ge 5$. On one hand, 
\begin{align*}
\binom{2p-5}{k+1}-\binom{p-2}{k+1}=\sum_{\ell=1}^{p-3}\binom{2p-5-\ell}{k} > (p-3)\binom{p-2}{k},
\end{align*}
and on the other hand, 
		\begin{align*}
	\binom{2p-3}{k}-\binom{p-1}{k+1} > \binom{p-2}{k}-\binom{p-1}{k+1}=-\binom{p-2}{k+1}=-\frac{p-2-k}{k+1}\binom{p-2}{k}.
		\end{align*}
Thus, we have  
\begin{align*}
    \f_k(d)& > (p-3)\binom{p-2}{k}-\frac{p-2-k}{k+1}\binom{p-2}{k}>0.
\end{align*}
	
(iv) We reason as in (i), and write
	\begin{equation*}
		\f_k(d):=\tau_k(2(d-1)+1,d-1)+\eta_{k-1}(2(d-1)+1,d-1)+\binom{d-1}{k}+\binom{d-2}{k}-\eta_k(2d+2,d).
	\end{equation*}
	We verify that $\f_{1}(4)=1,\f_{2}(4)=2$. We prove the statement by induction on $d\ge 4$ for all $1\le k \le d-2$. Assume that $d\ge 5$. Applications of  Lemma 9(iv) give that 
	\begin{align*}
		\f_k(d)=&\brac*{\binom{d}{k+1}+\binom{d-1}{k+1}+\binom{d-2}{k+1} -\binom{d-\floor{(d+1)/2}}{k+1}-\binom{d-\floor{(d+1)/2}-1}{k+1}}\\
		&+\brac*{\binom{d}{k}+2\binom{d-1}{k}-\binom{d-2}{k}}+\binom{d-1}{k}+\binom{d-2}{k}\\
		&-\left[\binom{d+1}{k+1}+2\binom{d}{k+1}-\binom{d-2}{k+1}\right]\\
		=&2\binom{d-2}{k+1}+2\binom{d-1}{k}
		-\binom{d}{k+1}-\binom{\ceil{(d-1)/2}}{k+1}-\binom{\ceil{(d-1)/2}-1}{k+1}.
	\end{align*}	
	The induction hypothesis on $d-1$ yields that $\f_k(d-1)>0$ and $\f_{k-1}(d-1)>0$ for all $1\le k\le d-3$.
	\begin{align*}
		\f_k(d-1)&=2\binom{d-3}{k+1}+2\binom{d-2}{k} -\binom{d-1}{k+1}-\binom{\ceil{d/2}-1}{k+1}-\binom{\ceil{d/2}-2}{k+1}\\
		\f_{k-1}(d-1)&=2\binom{d-3}{k}+2\binom{d-2}{k-1} -\binom{d-1}{k}-\binom{\ceil{d/2}-1}{k}-\binom{\ceil{d/2}-2}{k}.
	\end{align*}
	If $d$ is even, then $\f_k(d)=\f_{k-1}(d-1)+\f_{k}(d-1)$, which implies that $\f_k(d)>0$ for $1\le k\le d-3$. In the case of odd $d$,  $\f_k(d)\ge \f_{k-1}(d-1)+\f_{k}(d-1)$, which implies that $\f_k(d)>0$ for $1\le k\le d-3$. For $k=d-2$,
	\begin{align*}
		\f_{d-2}(d)=&2\binom{d-2}{d-1}+2\binom{d-1}{d-2}
		-\binom{d}{d-1}-\binom{\ceil{(d-1)/2}}{d-1}-\binom{\ceil{(d-1)/2}-1}{d-1}\\
		=&0+2(d-1)-d-0-0\\
		=&d-2>0.
	\end{align*}	
	This part is now proved.

(v) We reason as in (i), and write
	\begin{equation*}
		\f_k(d):=\tau_k(2(d-1)+1,d-1)+\tau_{k-1}(2(d-1)+1,d-1)+\binom{d-1}{k}+\binom{d-2}{k}-\eta_k(2d+2,d).
	\end{equation*}

Applications of  Lemma 9(iv) give that 
	\begin{align*}
		\f_k(d)=&\brac*{\binom{d}{k+1}+\binom{d-1}{k+1}+\binom{d-2}{k+1} -\binom{\ceil{(d-1)/2}}{k+1}-\binom{\ceil{(d-1)/2}-1}{k+1}}\\
		&+\brac*{\binom{d}{k}+\binom{d-1}{k}+\binom{d-2}{k} -\binom{\ceil{(d-1)/2}}{k}-\binom{\ceil{(d-1)/2}-1}{k}}\\
		&+\binom{d-1}{k}+\binom{d-2}{k}
		-\left[\binom{d+1}{k+1}+2\binom{d}{k+1}-\binom{d-2}{k+1}\right]\\
		=&\binom{d+1}{k+1}+\binom{d}{k+1}+\binom{d-1}{k+1} -\binom{\ceil{(d-1)/2}+1}{k+1}-\binom{\ceil{(d-1)/2}}{k+1}\\
		&+\binom{d-1}{k}+\binom{d-2}{k}
		-\binom{d+1}{k+1}-2\binom{d}{k+1}+\binom{d-2}{k+1}\\
		=&\binom{d-1}{k+1}-\binom{\ceil{(d-1)/2}+1}{k+1}-\binom{\ceil{(d-1)/2}}{k+1}.
	\end{align*}	
We show that  $\f_k(d)$ is strictly increasing in $d$ for all $k$: For $p\ge 2$, 
\begin{align*}
    \f_k(2p)&=\binom{2p-1}{k+1}-\binom{p+1}{k+1}-\binom{p}{k+1}\\
    \f_k(2p+1)&=\binom{2p}{k+1}-\binom{p+1}{k+1}-\binom{p}{k+1},
\end{align*}
which implies 
\begin{align*}
\f_k(2p+1)-\f_k(2p)&=\binom{2p}{k+1}-\binom{2p-1}{k+1}=\binom{2p}{k}>0\\
\f_k(2p+2)-\f_k(2p+1)&=\binom{2p+1}{k+1}-\binom{2p}{k+1}-\binom{p+2}{k+1}+\binom{p}{k+1}\\
&=\binom{2p}{k}-\binom{p+1}{k}-\binom{p}{k}>0.
\end{align*}
Since $\f_k(d)$ is increasing in $d$  for all $k$, and $\f_{1}(5)=2,\f_{2}(5)=3, \f_{3}(5)=1$, it follows that $\f_k(d)>0$ for all $k\in [1\ldots d-2]$.
	This part is now proved, which completes the proof of the lemma.
\end{proof}

\section{Acknowledgments}

Guillermo would like to express his gratitude to Julien Ugon for the stimulating discussions on the topic and for reviewing an earlier version of the paper.

\end{document}